\documentclass[12pt]{amsart}
\usepackage{comment}
\usepackage{macros}

\begin{document}
\title
[Resolution Property and the Brauer map] 
{The Resolution Property via Azumaya Algebras}

\author{Siddharth Mathur}
\address{Mathematisches Institut\\Heinrich-Heine-Universit\"at\\40204 D\"usseldorf, Germany.\\}
\email{siddharth.p.mathur@gmail.com} 

\begin{abstract}

\noindent Using formal-local methods, we prove that a separated and normal tame Artin surface has the resolution property. By proving that normal tame Artin stacks can be rigidified, we ultimately reduce our analysis to establishing the existence of Azumaya algebras. Our construction passes through the case of tame Artin gerbes, tame Artin curves, and algebraic space surfaces, each of which we establish independently. 

\end{abstract}
\maketitle

\section{Introduction}
\label{sec:intro}

An algebraic stack $X$ has the resolution property if every quasi-coherent sheaf of finite-type is the quotient of a vector bundle. It often happens that $X$ has enough line bundles (e.g. $X$ is quasi-projective or a smooth separated scheme) to resolve quasi-coherent sheaves. However, there are many important schemes and algebraic stacks that do not have enough line bundles, and for these spaces, the resolution property is very difficult to establish and poorly understood. The purpose of this paper is to prove the following

\begin{Theorem} \label{T1} Let $\ms X$ be a separated, tame Artin stack of finite-type over a field $k$.

\begin{enumerate}

\item If $\ms X$ is $1$-dimensional, then it satisfies the resolution property.

\item If $\ms X$ is $2$-dimensional and normal, then it satisfies the resolution property. \end{enumerate}

\noindent
In particular, they are global quotient stacks.

\end{Theorem}

In fact, our methods produce stronger results (see Theorems \ref{tamerigidification}, \ref{ResSurface}, and \ref{BrOrbisurface}). The strategy is to use the process of rigidification to reduce to the case of a stack with generically trivial stabilizers. In this setting, we are able to use formal-local methods to construct enough vector bundles. However, rigidification forces us to consider the case of gerbes. To deal with tame Artin gerbes we use non-abelian cohomology to further reduce the problem to finding Azumaya algebras. Then we generalize an argument of Gabber to prove:

\begin{Theorem} \label{T2} Let $\ms X$ be a separated, tame Artin stack which is finite-type over a field $k$. 
\begin{enumerate}

\item If $\ms X$ is $1$-dimensional, then $\Br(\ms X)=\Br'(\ms X)$.

\item If $\ms X$ is $2$-dimensional, regular in codimension $1$ ($R_1$), and is a gerbe over a stack with generically trivial stabilizers (e.g. $\ms X$ is a normal tame Artin surface) then $\Br(\ms X)=\Br'(\ms X)$.
\end{enumerate}
\end{Theorem}

Resolving quasi-coherent sheaves by vector bundles is a well-worn technique with a long history in algebraic geometry. Indeed, the following questions are some examples where it plays a crucial role.

\begin{enumerate}
\item (\cite[Section 2]{TotaroRes}) Does the $\text{K}$-theory of vector bundles agree with the $\text{K}$-theory of perfect complexes? 
\item (\cite{ConradDual}) Is every smooth, affine group scheme $G \to S$ linear? i.e. is there a closed immersion $G \to \text{GL}_S(V)$ for a vector bundle $V$ on $S$?
\item (\cite[Question 2]{Seshadri} and \cite[Theorem 3.8]{ThomasonRes}) Suppose $B$ is a ring of finite-type and $G$ a reductive group scheme over a general base scheme $S$. Is the ring of invariants, $B^G$, also of finite-type over $S$?
\item Is every cohomological Brauer class $\alpha \in \text{H}^2(X, \mathbf{G}_m)_{\text{tors}}$ representable by an Azumaya algebra?
\end{enumerate}
 
The resolution property can be viewed as a generalization of quasi-projectivity. Indeed, a line bundle on a quasi-compact scheme is ample if and only if the corresponding $\mathbf{G}_m$-bundle is quasi-affine. On the other hand, a striking characterization of the resolution property is the following:

\begin{Theorem} (Totaro, \cite{TotaroRes} and Gross, \cite{GrossRes}) \label{TotaroGross} Let $\ms X$ be a quasi-compact, quasi-separated (qcqs) algebraic stack with affine stabilizers at closed points. Then $\ms X$ satisfies the resolution property if and only if it can be expressed as a quotient stack $[W/\mathbf{GL}_n]$ where $W$ is quasiaffine. 
\end{Theorem}

Thus, the resolution property is closely related to another fundamental question: when is an algebraic stack \emph{globally} the quotient of a scheme? Although tame algebraic stacks are known to have this property locally, we cannot find a single separated Deligne-Mumford stack which is not a global quotient stack. In the instance of $\mu_n$- or $\mathbf{G}_m$-gerbes, this question is deeply intertwined with an important problem of Grothendieck.

In \cite{Brauer2}, Grothendieck posed the question: does every cohomological Brauer class arise from a $\text{PGL}_n$-torsor? In other words, is the Brauer map $\Br(X) \to \Br'(X)$ surjective? It is a difficult result due of Gabber that the Brauer map is surjective when $X$ admits an ample line bundle (see \cite{dejongample} for a proof). Using the resolution property we may reinterpret this result as follows: every $\mu_n$-gerbe morphism $X \to M$ satisfies the resolution property whenever $M$ admits an ample line bundle (see \cite{EHKV}, \cite{Gabbersthesis}, \cite{dejongample}). A solution to Grothendieck's question has remained out of reach for almost 50 years and it is not clear what one should expect, indeed see \cite{EHKV} for a non-separated counterexample. Note that the problem is completely open even when $M$ is a smooth scheme of dimension three, finite-type and separated over the complex numbers.

It is perhaps unexpected that the existence of Azumaya algebras should be tied to the resolution property of a very particular gerbe morphism (see Theorem \ref{BrauerMapTFAE}), but the story takes a further twist. Gabber's result can be used to verify the resolution property for a large class of Deligne-Mumford stacks that are \emph{not} gerbes. Indeed, Kresch and Vistoli \cite{KreschVistoli} proved that if a (generically tame) smooth Deligne-Mumford stack which is separated and finite-type over a field has quasi-projective coarse moduli space, then it is a global quotient stack. Their contribution is to give a method to reduce the general situation to the case of $\mu_n$-gerbes, whence Gabber's result applies via \cite{EHKV}. Their techniques are a mix of Giraud's nonabelian theory (see \cite{giraud}) and projective methods a la Bertini. One may conclude that such stacks have the resolution property by then resorting to \cite{TotaroRes} or \cite{GrossRes}. We give a criterion for smooth (or normal) tame Artin stacks to have the resolution property which is related to, but orthogonal, to their criterion (compare with \cite[Theorem 2]{KreschVistoli}).

\begin{Theorem} \label{TAG} Let $\ms X$ denote a tame Artin stack which is finite-type and separated over a field $k$. Moreover, suppose that $\ms X$ has a quasi-projective coarse moduli space.
\begin{enumerate} 
\item Let $\ms X$ be a gerbe, then $\ms X$ has the resolution property. 

\item The following are equivalent \begin{enumerate} \item Every smooth (resp. normal) $\ms X$ has the resolution property.
\item Every smooth (resp. normal) $\ms X$ which has generically trivial stabilizers has the resolution property. \end{enumerate} 
\end{enumerate}
\end{Theorem}

\noindent \textbf{Previous results in the literature}: Note that Theorem \ref{BrOrbisurface} and \ref{ResSurface} generalize the main results of Schr\"oer \cite{enoughAzumaya}, Schr\"oer-Vezzosi \cite{NormalSurfaceRes}, and Gross \cite{SurfaceRes}. Indeed, our results imply there are enough vector bundles on surfaces \emph{and} enough Azumaya algebras on $R_1$-surfaces. It is worth noting that in Schr\"oer-Vezzosi \cite{NormalSurfaceRes} and Gross \cite{SurfaceRes}, the authors do not address the case of algebraic spaces which are not schemes. Indeed, the existence of affine \emph{open} neighborhoods is crucial in several of their arguments and so new ideas are needed to tackle the algebraic space case. Thus, even if we assume $\ms X$ is an algebraic space, Theorem \ref{T1} is already new.

It should also be noted that the results of Kresch-Vistoli \cite{KreschVistoli} apply only when $\ms X$ is smooth, Deligne-Mumford and its coarse moduli space admits an ample line bundle. We have none of these hypothesis in Theorem \ref{T1} and this makes the problem much more difficult because we cannot leverage the existence of \emph{any} natural vector bundles (e.g. polarizing, cotangent or jet bundles) and in fact we must actively \emph{constuct} vector bundles instead. Very roughly, we accomplish this by glueing local vector bundles along analytic neighborhoods using the classical $\text{K}$-theory of commutative rings.
\\\\
\noindent \textbf{The structure of the paper}: In section one, we define the (relative) resolution property, the Brauer map, and show how these two topics are intertwined. There are some results here that do not appear in the literature but the expert can safely skip this section and return when necessary. The second section consists of a smorgasbord of auxiliary results that are required to prove Theorem \ref{T1}, we present them here because they are new and of independent interest. First, we give a deformation and obstruction theory for vector bundles on general algebraic stacks and then specialize to the tame setting where we prove results on the rigidity of vector bundles. Next, we introduce the notion of a Mayer-Vietoris square and use the theorem of Tannaka Duality to prove that flat neighborhoods often exist. Lastly, we prove that normal tame Artin stacks over a field can always be rigidified. In section three we prove the resolution property holds for algebraic space surfaces and tame Artin curves. Then, after a short study of tame Artin bands, we prove Theorem \ref{TAG}. In section four, we use formal glueing techniques to construct twisted vector bundles on certain gerbes and consequently prove Theorem \ref{T2}. Finally, we collect the results to establish Theorem \ref{T1}.
\\\\
\noindent \textbf{Acknowledgements} A large part of this paper comes from the authors Ph.D thesis. I would like to thank my advisor, Max Lieblich, for his generous support and guidance. I am also indebted to Jarod Alper and Aise Johan de Jong for their encouragement and suggestions. I also benefited from helpful discussions with Dan Abramovich, Dan Edidin, Lucas Braune, Jack Hall, Ariyan Javanpeykar, Andrew Kresch, Minseon Shin, David Rydh, Stefan Schr\"{o}er, Jason Starr, and Angelo Vistoli. I would also like to thank the referee for carefully reading an earlier version of this paper and providing many helpful comments. A part of this research was conducted in the framework of the research training group \emph{GRK 2240: Algebro-geometric Methods in Algebra, Arithmetic and Topology}, which is funded by the DFG.

\tableofcontents


\subsection{The Resolution Property and Quotient Stacks}

In this section we recall precise definitions of the (relative) resolution property (following Gross \cite{GrossRes} and \cite{Grossthesis}) and quotient stacks before explaining some basic consequences of the definitions. At the end, we relate the two concepts using results of Totaro and Gross.

\begin{definition} \label{Dgeneratingfamily}

Let $X$ be a quasi-compact and quasi-separated (qcqs) algebraic stack and $\mathscr{G}=\{F_i\}_{i \in I}$ a set of finitely-presented quasi-coherent sheaves. Then $\mathscr{G}$ is said to be a \emph{generating family} of $\mathcal{O}_X$-modules if for any quasi-coherent sheaf $M$ on $X$ there exists a surjection $\bigoplus F \to M$ where every $F \in \mathscr{G}$.

\end{definition}

A relative version of the above is

\begin{definition} \label{Df-generating} Let $f: X \to Y$ be a qcqs morphism of algebraic stacks. We say a family of finitely-presented quasi-coherent $\mathcal{O}_X$ modules $\mathcal{G}=\{F_j\}_{j \in J}$ is \emph{$f$-generating} if every quasi-coherent $\mathcal{O}_X$-module admits a surjection
\[\bigoplus_{(i,j) \in I \times J} F_j \otimes f^*N_i \to M\]
for some family of quasi-coherent $\mathcal{O}_Y$-modules $\{N_i\}_{i \in I}$ and where $F_j \in \mathscr{G}$.

The family is said to be \emph{universally $f$-generating} if $\{F_j|_{X \times_Y T} \}_{j \in J}$ is $f_T: X \times_Y T \to T$-generating for all morphisms $T \to Y$.
\end{definition}

Now we can define the resolution property

\begin{definition} \label{Dresolutionproperty}

Let $f: X \to Y$ be a qcqs morphism of algebraic stacks. We say $f$ has the \emph{resolution property} if there exists a universally $f$-generating family of finitely-presented locally free $\mathcal{O}_X$-modules $\{V_i \}_{i \in I}$. We say that $X$ has the \emph{resolution property} if there exists a generating family of finitely-presented locally free $\mathcal{O}_X$-modules $\{V_i\}_{i \in I}$.

\end{definition}

\begin{remark} \label{RemarkoverZ} For a qcqs algebraic stack $X$, the resolution property of the canonical morphism $f: X \to \Spec \mathbf{Z}$ is equivalent to the resolution property of $X$. Indeed, if a universally $f$-generating family of locally free sheaves on $X$ exists, then every quasi-coherent module $F$ can be surjected onto by a direct sum of vector bundles $\bigoplus V_i^{\oplus m_i} \to F$. Conversely if $X$ satisfies the resolution property then take a generating family of locally free sheaves $\{V_i\}$, it is certainly $f$-generating hence it is \emph{universally} $f$-generating by \cite[Corollary 2.10]{GrossRes} and the fact that $\Spec \mathbf{Z}$ has affine diagonal.
\end{remark}

The relativization will come in handy because if a qcqs algebraic stack $X$ maps to $Y$, to show that $X$ satisfies the resolution property it suffices to show that $X \to Y$ and $Y$ both satisfy the resolution property: 

\begin{lemma} \label{resolutionpropertycomposes} Suppose $f: X \to Y$ is a qcqs morphism between qcqs algebraic stacks. Then $X$ satisfies the resolution property if $f$ and $Y$ both satisfy the resolution property.
\end{lemma}

The proof is standard so we omit it, a more general result is given in \cite[Proposition 2.8 (v)]{GrossRes}. One of the difficulties of the resolution property is that we do not know when it descends. However, the following is known.



\begin{lemma} (\cite{GrossRes} Proposition 5.3 (vii)) \label{finiteflatdescent} Suppose $f: X \to Y$ is a finite, faithfully flat and finitely-presented morphism between qcqs $S$-stacks. Then if $X \to S$ satisfies the resolution property, $Y \to S$ does as well.
\end{lemma}

Now we introduce the notion of a quotient stack following \cite{EHKV}.

\begin{definition} \label{Dquotientstack} Let $X$ be a qcqs algebraic stack, we say that $X$ is a \emph{quotient stack} if it is isomorphic to $[Z/\mathbf{GL}_n]$ where $Z$ is an algebraic space.
\end{definition}

In general, it is very difficult to verify if an arbitrary algebraic stack is a quotient stack but we give a criterion below. Recall that a morphism of algebraic stacks $X \to Y$ is said to be projective if there exists a finite-type quasi-coherent sheaf $E$ on $Y$ and a factorization $X \to \mathbf{P}(E)_Y \to Y$ where $\mathbf{P}(E)_Y$ is the projectivization of $E$ over $Y$ and $X \to \mathbf{P}(E)_Y$ is a closed immersion. Also, let $I_X \to X$ denote the inertia stack and for any morphism $f: T \to X$, let $I_f \to T$ denote the associated base change.

\begin{proposition} (\cite[Lemmas 2.12 and 2.13]{EHKV}) \label{quotientstackTFAE} Let $X$ be an algebraic stack finite-type over a noetherian base scheme $S$, then the following are equivalent
\begin{enumerate}

\item $X$ is a quotient stack

\item There exists a locally free sheaf of finite rank $V$ on $X$ so that for every geometric point $x: \Spec k \to X$ the stabilizer action $I_x$ on the vector space $V_{\Spec k, x}$ is faithful i.e. the morphism of group schemes $I_x \to \mathbf{GL}(V_{\Spec k, x})$ is injective.

\item There exists a faithfully flat, projective morphism $Y \to X$ where $Y$ is a quotient stack.

\end{enumerate}

\end{proposition}

\begin{remark} \label{RemarkLinearArbitrary} At first glance, it may seem very restrictive to only consider stacks of the form $[Z/\mathbf{GL}_n]$ but in fact this class of stacks includes all those of the form $[Y/G]$ where $Y$ is an algebraic space and $G$ a flat linear algebraic group (i.e those embeddable into $\mathbf{GL}_n$). Indeed, because $G$ is linear, there is a closed embedding $G \subset \mathbf{GL}_n$ and we may consider the contracted product $Y \times^{G} \mathbf{GL}_n=[Y \times \mathbf{GL}_n/G]$ and view it as a $\mathbf{GL}_n$-space via right translation. Moreover $[Y/G] \cong [Y \times^{G} \mathbf{GL}_n/\mathbf{GL}_n]$ (see \cite[Remark 2.11]{EHKV}).
\end{remark}

\begin{definition} \label{Dfaithfulvectorbundle} Let $\mathscr{X}$ be a qcqs algebraic stack and $V$ a vector bundle on $\ms X$. If for every geometric point $x: \Spec k \to \ms X$ the stabilizer action $I_x$ on the vector space $V_{\Spec k, x}$ is faithful, then we say $V$ is a \emph{faithful vector bundle}. Let $\ms X \to \ms Y$ denote a qcqs morphism of algebraic stacks and $V$ a vector bundle on $\ms X$. Suppose $V$ is faithful after restriction along a smooth schematic cover $Y \to \ms Y$, then $V$ is said to be a \emph{relatively faithful vector bundle}. 
\end{definition}

It is not difficult to prove the following:

\begin{proposition} \label{faithfulvectorbundleTFAE} Let $\mathscr{X}$ be a qcqs algebraic stack and let $V$ be a vector bundle. Then $V$ is a faithful vector bundle if and only if $\text{Frame}(V)$ is an algebraic space.
\end{proposition}

Theorem \ref{TotaroGross} relates these two seemingly disparate notions: a qcqs algebraic stack with affine stabilizers at closed points satisfies the resolution property if and only if it admits a quotient stack presentation of the form $[W/\mathbf{GL}_n]$ where $W$ is a quasi-affine scheme. Another relationship between these two notions is the following:

\begin{proposition} \label{resolutionpropertycoarsemoduli} Let $X$ be an algebraic stack which is finite-type over a noetherian base scheme $S$. Assume moreover that $X$ has finite diagonal. If $\pi: X \to M$ is the coarse moduli space of $X$, then $X$ is a quotient stack if and only if $\pi$ has the resolution property.
\end{proposition}

\begin{proof} Suppose that $X=[Z/\mathbf{GL}_n]$ for an algebraic space $Z$. We will show that $Z \to M$ is an affine morphism then use \cite[Theorem 6.10]{GrossRes} to conclude that $\pi$ has the resolution property. Note that by \cite[Theorem 2.7]{EHKV}  there exists a finite surjective morphism $f: Y \to X$ where $Y$ is a scheme. It follows that $\pi \circ f: Y \to M$ is also finite. Suppose that $M$ is affine, then $Y$ is affine and because $Z \times_X Y$ is a $\mathbf{GL}_n$-torsor over $Y$, it too must be affine. However, $Z \times_X Y \to Z$ is a finite surjective morphism so by Chevalley's theorem (see \cite[Tag 07VP]{stacks}) $Z$ must be affine as well.

Conversely, if $\pi$ has the resolution property, \cite[Theorem 6.10]{GrossRes} implies there is a quasiaffine classifying map $X \to \B\mathbf{GL}_{n,M}$ over $M$. Therefore, pulling back the universal frame bundle yields a $\mathbf{GL}_n$-torsor $Z \to X$ where $Z$ is quasiaffine over $M$. In particular, $Z$ is an algebraic space and therefore $X=[Z/\mathbf{GL}_n]$ is a quotient stack. \end{proof}

\subsection{Grothendieck's Question on the Brauer map}
We begin by briefly reviewing the theory of the Brauer group with a focus on its relationship to gerbes and twisted sheaves. For what follows, let $X$ denote an algebraic stack endowed with the fppf topology.

\begin{definition} \label{DefCohoBrauer} The \emph{cohomological Brauer group} of $X$ is $\H^2(X,\mathbf{G}_m)_{\text{tors}}$, we denote it by $\Br'(X)$.
\end{definition}

Note that the following sequence is exact in the category of sheaves groups on $X$ (see \cite[pg. 341]{giraud}).
\[1 \to \mathbf{G}_m \to \mathbf{GL}_n \to \mathbf{PGL}_n \to 1\]
Following Giraud's nonabelian formalism \cite{giraud}, we obtain an exact sequence of pointed sets for any $n$:
\[ ... \to \H^1(X,\mathbf{G}_m) \to \H^1(X,\mathbf{GL}_n) \to \H^1(X,\mathbf{PGL}_n) \to \H^2(X,\mathbf{G}_m)\]
Here, $\H^1(X, G)$ denotes isomorphism classes of left $G$-torsors and $\H^2(X, \mathbf{G}_m)$ denotes isomorphism classes of $\mathbf{G}_m$-gerbes. Let $\delta_n$ denote the connecting map from degree one to two.
\begin{definition} \label{DefBrauer} The \emph{geometric Brauer group} is defined to be $\bigcup_{n \geq 1} \text{Im}(\delta_n)$, we denote it by $\Br(X)$.
\end{definition}
Noting that $\mathbf{PGL}_n$-torsors correspond to Azumaya algebras of degree $n$ we obtain

\begin{Theorem} \label{BrauerMap} (\cite[Chapter V.4.4]{giraud}) Let $X$ be an algebraic stack, then $\Br(X)$ is a group and there is a natural injection $\delta: \Br(X) \to \Br'(X)$, called the \emph{Brauer map}. The map sends the class of an Azumaya Algebra $A$ to the $\mathbf{G}_m$-gerbe of trivializations of $A$, $\mathscr{X}_A$. More precisely, $\mathscr{X}_A(T)$ is the category of pairs $(V, \phi)$ where $V$ is a vector bundle on $T$ and a trivialization $\phi: \End(V) \cong A|_T$. A morphism from $(V,\phi)$ to $(W, \psi)$ is an isomorphism $V \to W$ compatible with the trivializations $\phi$ and $\psi$.
\end{Theorem}

Now we can formulate a problem posed by Grothendieck in \cite{Brauer2}:

\begin{Question} \label{GrothendieckQuestion} (Grothendieck) Is the Brauer map $\delta: \Br(X) \to \Br'(X)$ surjective?
\end{Question}

\noindent In other words, does every torsion cohomology class come from an Azumaya algebra? 

In what follows, we define twisted sheaves which, at the expense of introducing stacky complexity, will clarify the question above. However, before doing so we briefly recall the notion of a $\emph{Band}$ (or \emph{Lien} in French) on a site $S$. Begin with the category of group sheaves on $S$ and note that for any two objects $G$ and $H$, the sheaf $\text{Isom}(G,H)$ admits an action by the inner automorphism group of $H$. Let $\ms B_S$ denote the category whose objects are groups sheaves on $S$ and whose morphisms are given by the values of the sheaf $\text{Isom}(G,H)/\text{Inn}(H)$, the category of bands on $S$ is then defined to be the stackification of $\ms B_S$. Note that given a group sheaf $G$, it has an associated band which we sometimes denote by $\text{Band}(G)$. There is a canonical way to associate a band to a gerbe $\ms X$ and we denote it by $\text{Band}(\ms X)$, moreover, if there is an isomorphism $\text{Band}(\ms X) \cong \text{Band}(G)$, then we say that $\ms X$ is \emph{banded} by $G$. For the definition of a band associated to a gerbe we refer the reader to \cite[Chapter IV, Section 1]{giraud} (or \cite[3.2,3.3]{EHKV} for a brief discussion). 

Now, suppose that $D \subset \mathbf{G}_m$ is a subgroup scheme, then the character group $\hat{D}$ is isomorphic to $\mathbf{Z}$ or $\mathbf{Z}/n\mathbf{Z}$ and any quasi-coherent sheaf $F$ on a gerbe banded by $D$ admits an eigendecomposition $F=\bigoplus_{\chi \in \hat{D}} F_{\chi}$ by \cite[Proposition 3.1.1.4]{lieblich2008twisted}.

\begin{definition} \label{DefTwisted} Let $\mathscr{X}$ be a $D$-gerbe on $X$. If $m \in \mathbf{Z}$ and $\bar{m} \in \hat{D}$ is its image, then a quasi-coherent sheaf $F$ is \emph{m-twisted} if its eigendecomposition is $F=F_{\bar{m}}$. If $F$ is a 1-twisted quasi-coherent sheaf, then we will just say that $F$ is a \emph{twisted} sheaf and if it is $0$-twisted we will say it is \emph{untwisted}.
\end{definition}

\begin{remark} \label{twistedfaithful} One may check that if $D=\mu_n$, then an $m$-twisted sheaf is faithful exactly when $(m,n)=1$. Moreover, if $D=\mathbf{G}_m$ then an $m$-twisted sheaf is faithful only when $m=1,-1$. Also note that the tensor product of a $n$-twisted sheaf and a $m$-twisted sheaf is $n+m$-twisted. \end{remark}

The following theorem explains the connection between the Brauer map, twisted sheaves and the resolution property.

\begin{Theorem} \label{BrauerMapTFAE}  Let $X$ be an quasicompact algebraic stack and $[\mathscr{X}] \in \H^2(X, \mathbf{G}_m)$ a $\mathbf{G}_m$-gerbe. Then the following statements are equivalent.
\begin{enumerate}
\item $\mathscr{X}$ is isomorphic to $\mathscr{X}_A$ for some Azumaya algebra $A$ on $X$ i.e. $[\mathscr{X}]$ is in the image of the Brauer map.
\item There exists a twisted vector bundle of nonzero constant rank on $\mathscr{X}$.
\item There exists a $\mu_n$-gerbe $\mathscr{Y}$ admitting a twisted vector bundle and a map of stacks $\mathscr{Y} \to \mathscr{X}$ equivariant for the inclusion $\mu_n \subset \mathbf{G}_m$. Moreover, any $\mu_n$-gerbe $\mathscr{Y}$ with a map $\mathscr{Y} \to \mathscr{X}$ equivariant for $\mu_n \subset \mathbf{G}_m$ admits a twisted vector bundle.
\item There exists a $\mu_n$-gerbe $\mathscr{Y}$ so that the structure map $\mathscr{Y} \to X$ has the resolution property and a map of stacks $\mathscr{Y} \to \mathscr{X}$ equivariant for the inclusion $\mu_n \subset \mathbf{G}_m$.
\item There is a finite, finitely presented flat map $Y \to X$ so that $\mathscr{X}|_Y $ is isomorphic to $\mathscr{X}_A$ for an Azumaya algebra $A$ on $Y$.
\end{enumerate}
In particular, to show the Brauer map $\Br(X) \to \Br'(X)$ is surjective, it suffices to show that all $\mu_n$-gerbe morphisms $\ms Y \to X$ have the resolution property.
\end{Theorem}
\begin{proof}

\underline{(1) $\Rightarrow$ (2)}

Let $T$ be any algebraic space over $X$, then $\mathscr{X}(T)$ is the groupoid of pairs $(V, \phi)$ where $V$ is a vector bundle on $T$ and $\phi$ is an isomorphism $\End(V) \to A|_T$. There is a tautological vector bundle on this stack: it is the sheaf $(V,\phi) \mapsto \Gamma(T, V)$ where the restriction maps are those induced by morphisms in $\mathscr{X}_A$. Noting that the automorphisms of $(V, \phi)$ are precisely left multiplication on $V$ by an element in $\mathbf{G}_m(T)$, (2) follows.

\underline{(2) $\Rightarrow$ (1)}

Let $V$ be a twisted vector bundle of nonzero constant rank and consider $\text{End}(V)$. Since the relative stabilizers of $\mathscr{X} \to X$ act trivially on this algebra, $\text{End}(V)$ is actually (the pullback of) an algebra on $X$. In fact, it is an Azumaya algebra on $X$. Indeed, \'etale locally on $X$, the gerbe map $\mathscr{X} \to X$ admits a section which, in turn, yields an equivalence $\mathscr{X} \to \B\mathbf{G}_m$. This equivalence endows $\mathscr{X}$ with a twisted line bundle, $L$, which in turn, determines a natural trivialization $\End(V) \cong \End(V \otimes L^{\vee})$ (using the natural identifications $W \otimes W^{\vee}=\End(W)$ and $L \otimes L^{\vee}=\mathcal{O}_{\ms X}$ for any vector bundle $W$ and line bundle $L$). Thus, any $T$-point of $\mathscr{X}$ determines a trivialization of $\End(V)|_T$. In other words, we obtain a map $\mathscr{X} \to \mathscr{X}_{\End(V)}$ which respects the $\mathbf{G}_m$-gerbe structure. Since this must be an isomorphism, the first statement now follows.

\underline{(2) $\Rightarrow$ (3)}

The previous argument shows that (2) implies $\mathscr{X}_A \cong \mathscr{X}$ for some Azumaya algebra $A$ on $X$. It follows that the corresponding cohomology class is torsion and hence admits a lift to a cohomological class in $\text{H}^2(X,\mu_n)$ via the Kummer sequence. By the nonabelian formalism of Giraud, this implies there exists a $\mu_n$-gerbe $\mathscr{Y}$ and a map $\mathscr{Y} \to \mathscr{X}$ equivariant for $\mu_n \to \mathbf{G}_m$. If $\mathscr{Y}$ is any $\mu_n$-gerbe admitting such a morphism, we may pull back a given twisted vector bundle on $\mathscr{X}$ to obtain one over $\mathscr{Y}$.

\underline{(3) $\Rightarrow$ (1)}

Suppose (3) holds and let $V$ denote a twisted vector bundle on the $\mu_n$-gerbe $\sY$. Set $A=\End(V)$ and observe that this is the pullback of an Azumaya algebra on $X$. Consider the $\mathbf{G}_m$-gerbe associated to it, $\sX_A$, whose $T$-points are $(W, \phi: \End(W) \cong A|_T)$. Observe that there is a natural map of gerbes $\sY \to \sX_A$ equivariant for $\mu_n \subset \mathbf{G}_m$. Indeed, a $T$-point of $\sY$ induces a trivialization $\sY_T \cong \B\mu_{n,T}$ i.e. a $T$-point of $\sY$ yields a twisted line bundle $L$ on $\sY_T$. But then $V|_{\sY_T} \otimes L^{\vee}$ is a vector bundle on $T$ and there is a natural isomorphism $\phi: \End(V|_{\sY_T} \otimes L^{\vee}) \cong \End(V)|_T=A|_T$, i.e. a trivialization of the Azumaya algebra $A$. In summary, a $T$-point of $\sY$ naturally induces a $T$-point of $\sX_A$. One can check the maps on stabilizer groups respect the inclusion $\mu_n \subset \mathbf{G}_m$. However, since there is also a map $\sY \to \sX$ equivariant for $\mu_n \subset \mathbf{G}_m$, this forces $\sX_A \cong \sX$.

\underline{(3) $\Rightarrow$ (4)}

Suppose (3) holds. Fix a twisted vector bundle $V$ on $\mathscr{Y}$, we will show that the associated frame bundle is quasi-affine over $X$. This implies $\mathscr{Y} \to X$ factors through a quasiaffine morphism $\mathscr{Y} \to \B\mathbf{GL}_{n,X}$ so by \cite[Theorem 6.10]{GrossRes}, it would follow that $\mathscr{Y} \to X$ satisfies the resolution property. Let $F$ denote the total space of the frame bundle of $V$ and let $F \to \mathscr{Y} \to X$ be the associated maps. Since $V$ is faithful, $F$ is an algebraic space and by passing to a smooth cover of $X$ we may suppose $X$ is a scheme and that $\mathscr{Y} \cong \B\mu_n$. Now $F$ is cohomologically affine over $X$. Indeed, $F \to \mathscr{Y}$ is a $\mathbf{GL}_n$-torsor and $\mathscr{Y}=\B\mu_n \to X$ is cohomologically affine and therefore $F \to X$ is as well. By \cite[Proposition 3.3]{goodmodulispaces}, this implies the morphism $F \to X$ is affine, as desired.

\underline{(4) $\Rightarrow$ (3)}

Now suppose (4) holds. Since $\mathscr{Y} \to X$ is a $\mu_n$-gerbe the relative inertia $I_f \to \mathscr{Y}$ has geometric fibers isomorphic to $\mu_n$. Thus the assumption that $f$ has the resolution property and \cite[Theorem 6.10]{GrossRes} implies there is a vector bundle $V$ on $\mathscr{Y}$ whose frame bundle is quasi-affine over $X$. We may decompose $V=\bigoplus_{i \in \mathbb{Z}/n\mathbb{Z}} V_i$ into its eigensheaves with respect to the action of $\mu_n$. Let $S \subset \{1,...,n-1\}$ be the subset consisting of integers where $V_i \neq 0$. If there was a $i \in S$ relatively prime to $n$ then there is a $k$ such that $V_i^{\otimes k}$ is $1$-twisted. Suppose every $i$ is not relatively prime to $n$. There cannot be a prime $p$ dividing every integer in $S$ and $n$, otherwise $\mu_p \subset \mu_n$ would act trivially on the fibers of $V$. This is contrary to the fact that $\mu_n$ acts faithfully on the geometric fibers of $V$. It follows that there exists a polynomial combination of the $V_i$'s that is $1$-twisted.

\underline{(1) $\Leftrightarrow$ (5)}

For the proof that 5. implies 1, we refer the reader to \cite[Chapter 2, Lemma 4]{Gabbersthesis}.  For the converse, set $Y=X$ and let $Y \to X$ be the identity. \end{proof}

\section{Tools}

\subsection{Deformation Theory of Vector Bundles on Tame Stacks}

We begin by explaining a deformation-obstruction theory of vector bundles on algebraic stacks and then we discuss the rigidity of representation types on Tame Artin stacks.

Fix a square zero $f: \mathscr{X}_0 \to \mathscr{X}$ thickening of algebraic stacks, that is, a closed immersion whose defining ideal is square zero. We will work over the lisse-\'etale site of $\ms X$.

\begin{definition} \label{Ddeformation} Let $V_0$ be a vector bundle on $\mathscr{X}_0$, a \emph{deformation} of $V_0$ to $\sX$ is a pair $(V, \phi)$ where $V$ vector bundle on $\sX$ and $\phi: V|_{\sX_0} \to V_0$ is an isomorphism.
\end{definition}

\noindent In \cite[Tag 08VW]{stacks} one can find a description of the deformations and obstructions of a module under a thickening of ringed topoi. This result does not apply to our current situation as a square zero thickening of algebraic stacks may not induce a thickening of topoi. Indeed, it is not clear that a square zero thickening of algebraic stacks induces an equivalence on lisse-\'etale topoi. One may circumvent this issue by using the results of Olsson (see \cite[Theorem 1.5]{olssondef}) to develop a deformation theory of \emph{faithful} vector bundles. Instead, we describe an analogous description of the deformation and obstructions of an arbitrary vector bundle under a square zero thickening of algebraic stacks.

\begin{proposition} \label{deformationobstruction} Let $f: \sX_0 \to \sX$ be a square zero thickening of algebraic stacks with defining ideal $I \subset \mathcal{O}_{\mathscr{X}}$ and $V_0$ a vector bundle on $\mathscr{X}_0$.
\begin{enumerate}

\item There exists an obstruction $o \in \text{Ext}^2_{\mathcal{O}_{\mathscr{X}_0}}(V_0, I \otimes V_0)$ whose vanishing is equivalent to the existence of a deformation of $V_0$.

\item If a deformation of $V_0$ exists, the set of isomorphism classes of pairs $(V, \phi)$ lifting $V_0$ is a (trivial) torsor under $\text{Ext}^1_{\mathcal{O}_{\mathscr{X}_0}}(V_0, I \otimes V_0)$.

\end{enumerate}

\end{proposition}

\begin{proof} This argument is standard so we sketch it briefly. One shows that the stack of deformations of $V_0$, $\mathscr{D}$, is a gerbe over $\sX_{\text{lis-\'et}}$ banded by the quasi-coherent sheaf $f_*\text{\underline{Hom}}_{\mathcal{O}_{\mathscr{\sX}_0}}(V_0, V_0 \otimes I)$. This suffices because equivalence classes of gerbes on $\sX_{\text{lis-\'et}}$ banded by an abelian sheaf $A$ correspond to cohomology classes in $\H^2(\sX_{\text{lis-\'et}}, A)$ (see \cite[Theorem 12.2.8]{Olsson} or \cite[Chapter 4, Theorem 3.4.2]{giraud}) and
\[\H^2(\sX_{\text{lis-\'et}}, f_*\text{\underline{Hom}}_{\mathcal{O}_{\mathscr{\sX}_0}}(V_0, V_0 \otimes I)) \cong \H^2(\sX_{0, \text{lis-\'et}}, \text{\underline{Hom}}_{\mathcal{O}_{\mathscr{\sX}_0}}(V_0, V_0 \otimes I)) \cong \text{Ext}^2_{\mathcal{O}_{\mathscr{X}_0}}(V_0, V_0 \otimes I)\]
Moreover, under this correspondence, a gerbe is trivial over $\sX_{\text{lis-\'et}}$ if and only if it corresponds to the trivial cohomology class. The first result now follows since a gerbe is trivial if and only if there is an equivalence $\mathscr{D} \cong \text{B}A_{\sX_{\text{lis-\'et}}}$ and such equivalences correspond to sections of $\mathscr{D}$ over $\sX_{\text{lis-\'et}}$ i.e. deformations of $V_0$. The second result follows because isomorphism classes of sections of $\mathscr{D}$ correspond to $f_*\text{Hom}_{\mathcal{O}_{\mathscr{\sX}_0}}(V_0, V_0 \otimes I)$-torsors which in turn correspond to cohomology classes in
\[\H^1(\sX_{\text{lis-\'et}}, f_*\text{\underline{Hom}}_{\mathcal{O}_{\mathscr{\sX}_0}}(V_0, V_0 \otimes I)) \cong \H^1(\sX_0, \text{\underline{Hom}}_{\mathcal{O}_{\mathscr{\sX}_0}}(V_0, V_0 \otimes I)) \cong \text{Ext}^1_{\mathcal{O}_{\mathscr{X}_0}}(V_0, V_0 \otimes I)\]

The category $\mathscr{D}$ is defined over $\sX_{\text{lis-\'et}}$: its objects over $T \in \sX_{\text{lis-\'et}}$ are deformations $(V, \phi)$ of $V_0|_{T_0}$ where $T_0=T \times_{\mathscr{X}} \mathscr{X}_0$. An arrow $a: (V, \phi) \to (V', \phi')$ over $T$ is an isomorphism $a: V \cong V'$ which is compatible with the isomorphisms $\phi$, $\phi'$. The pullback of an object $(V, \phi) \in \mathscr{D}(T)$ along a morphism $t: T' \to T$ in $\sX_{\text{lis-\'et}}$ is simply $(t^*V, t^*\phi)$. The proof that $\mathscr{D}$ is a gerbe $\sX_{\text{lis-\'et}}$ banded by $f_*\text{\underline{Hom}}_{\mathcal{O}_{\mathscr{\sX}_0}}(V_0, V_0 \otimes I)$ is standard so we omit it. \end{proof}

Below we begin a discussion of deformation-theoretic phenomena of vector bundles on tame Artin stacks which are gerbes (or, more succinctly, tame Artin gerbes).

\begin{lemma} \label{vectorbundlesrigid} Consider a tame Artin gerbe $\mathscr{X} \to \spec A$ over a local Artin ring with residue field $k$. Let $V$ and $W$ be vector bundles on $\mathscr{X}$. If there is an isomorphism 
\[V|_{\mathscr{X} \times_{\Spec A} \Spec K} \cong W|_{\mathscr{X} \times_{\Spec A} \Spec K}\]
for some field extension $k \subset K$, then there exists an isomorphism $V \cong W$ on $\mathscr{X}$.
\end{lemma}

\begin{proof}

Consider the sheaves $\text{\underline{Isom}}_{\mathscr{X}}(V, W)=I \subset \text{\underline{Hom}}_{\mathscr{X}}(V, W)=H$ over $\Spec A$. Since gerbes are always flat and locally of finite presentation over their coarse space (see \cite[Tag 06QI]{stacks}) it follows that $H$ is representable by a finitely-presented abelian cone over $\Spec A$, (see \cite[Theorem D]{Hallcoho}). Our hypothesis is that the open subfunctor, $I$, admits a $\Spec K$-point.

Since $H \times_{\Spec A} \Spec k$ is a abelian cone of finite-type over a field, it must be isomorphic to $\mathbf{A}^n_k$. By hypothesis the open subfunctor $I \times_{\Spec A} \Spec k \subset H \times_{\Spec A} \Spec k \cong \mathbf{A}^n_k$ is nonempty after extending coefficients to the larger field $K$. This implies $I \times_{\Spec A} \Spec k$ is nonempty. However, nonempty open subschemes of $\mathbf{A}^n_k$ always have a $k$-rational point when $k$ is infinite i.e. there exists an isomorphism $\phi: V|_{\mathscr{X} \times_{\Spec A} \Spec k} \cong W|_{\mathscr{X} \times_{\Spec A} \Spec k}$. If $k$ is finite then we exploit Lang's theorem: every torsor for a smooth
connected algebraic group over a finite field is trivial (see \cite[Theorem 2]{lang1956algebraic}). Note that $I$ is a torsor under $\text{Aut}(W)$. Moreover, $\text{Aut}(W)$ is open in the \emph{affine} space $\text{\underline{Hom}}_{\sX}(W,W)$ and therefore smooth and connected. It follows that $I$ admits a $k$-point even if $k$ is finite.

Thus, $(W, \phi)$ is a deformation of $V|_{\mathscr{X} \times_{\Spec A} \Spec k}$ to $\mathscr{X}$. Since $(V, \text{can})$ is another deformation (here "can" denotes the canonical isomorphism), the lemma will follow once we show deformations are unique.

Factor the map $A \to k$ into square zero extensions so that we may use Proposition \ref{deformationobstruction}. Let $J$ be an ideal sheaf defining a square-zero thickening. To show the deformation spaces vanish:
\[\H^1(\mathscr{X} \times_{\Spec A} \Spec k, J \otimes \text{End}(V))=0\]
observe that the Leray spectral sequence for the coarse moduli space map $\pi: \mathscr{X} \times_{\Spec A} \Spec k \to \Spec k$ degenerates by tameness. Indeed, tameness implies the higher direct images of $J \otimes \End(V)$ vanish so the Leray spectral sequence yields an isomorphism: 
\[\H^1(\mathscr{X} \times_{\Spec A} \Spec k, J \otimes \text{End}(V))=\H^1(\Spec k, \pi_*(J \otimes \End(V)))\] 
Thus, the deformation spaces vanish identically since they correspond to the cohomology of a coherent sheaf over a point. The result follows. \end{proof}

\begin{lemma} \label{vectorbundlesrigidlocalring} Let $\pi: \mathscr{X} \to \Spec R$ be a coarse space morphism from a tame Artin stack to a local Noetherian ring with residue field $k$. Fix two vector bundles on $V$ and $W$ on $\mathscr{X}$ with an isomorphism over the closed fiber $V_k \cong W_k$, then $V \cong W$ over $\ms X$. 
\end{lemma}

\begin{proof} Denote the closed fiber by $\mathscr{X} \times_X \Spec k = \sX_k$ and the closed immersion by $\iota: \sX_k \to \sX$. We will show that the natural adjunction morphism of coherent sheaves
\[\underline{\text{Hom}}_{\mathcal{O}_{\sX}}(V,W) \to \iota_* \underline{\text{Hom}}_{\mathcal{O}_{\sX_k}}(V_k, W_k) \]
is surjective on global sections. It is surjective as a morphism of sheaves. Note that there is an equality of functors $\Gamma(\sX, -) = \Gamma(\Spec R, -) \circ \pi_*(-)$, and both $\pi_*$ and $\Gamma(\Spec R, -)$ are exact functors on $\text{Qcoh}(\sX)$ and $\text{Qcoh}(\Spec R)$ respectively. It follows that $\Gamma(\sX, -)$ sends a surjection to a surjection, so that
\[\text{Hom}_{\mathcal{O}_{\sX}}(V,W) \to \iota_* \text{Hom}_{\mathcal{O}_{\sX_k}}(V_k, W_k)\]
is surjective. Thus, we may lift the isomorphism $V_k \cong W_k$ to a homomorphism $\phi: V \to W$. Since it is a isomorphism on the closed fiber, Nakayama's lemma implies it is surjective. However, since they are vector bundles of the same rank, the surjection $\phi$ must be an isomorphism. \end{proof}

\noindent The following proposition is straightforward to show:

\begin{proposition} \label{twistedopenclosed} Let $\sX \to X$ be a gerbe banded by $\mu_n$ and $V$ a vector bundle on $\sX$. The locus where $V$ is a twisted vector bundle is open and closed.
\end{proposition}


\subsection{Mayer Vietoris squares and Flat Neighborhoods via Tannakian Duality}

\begin{definition} \label{DMVsquares} Consider a cartesian square of algebraic stacks, where $i$ is an open immersion:

\[\begin{tikzcd}
U \times_X Y \rar \dar & Y \dar{f} \\
U \rar{i} & X \\
\end{tikzcd}\]

\noindent In addition, suppose that $f$ is flat and $Z \times_X Y \to Z$ is an isomorphism for a finitely presented closed immersion $Z \to X$ where $X \backslash Z=U$. Such a square is called a \emph{flat Mayer-Vietoris square} and $Y$ is called a \emph{flat neighborhood} of $Z$. 
\end{definition}

\begin{remark} In \cite{hall2016mayer} the authors present a definition for a flat Mayer-Vietoris square without any finiteness hypothesis. However, in the absence of this finiteness condition, the authors require an additional hypothesis: they ask that for any $W \to X$ whose image is disjoint from $U$, the map $W \times_X Y \to W$ is an isomorphism. However, \cite[Lemma 3.2(2)]{hall2016mayer} guarantees that this is equivalent to our definition when $Z$ is finitely-presented.  \end{remark}

The property of being a flat Mayer-Vietoris square is flat local on $X$:

\begin{lemma} (\cite[Lemma 3.1]{hall2016mayer}) \label{MVsquaresflatlocal} Suppose we have a flat Mayer-Vietoris square and a map of algebraic stacks $X' \to X$. We may base change the square to $X'$ to obtain

\[\begin{tikzcd}
U' \times_{X'} Y' \rar \dar & Y' \dar{f'} \\
U' \rar{i'} & X' \\
\end{tikzcd}\]

\noindent then this new square is also a flat Mayer-Vietoris square. Moreover, the property of being a flat Mayer-Vietoris square is flat-local on $X$. \end{lemma}

The utility of flat Mayer-Vietoris squares is that they give us a new way of constructing vector bundles on schemes and stacks. Indeed, consider the following glueing result:

\begin{Theorem} \label{MVglueing} (\cite[Theorem B(1)]{hall2016mayer}) Suppose we have a flat Mayer-Vietoris square as in Definition \ref{DMVsquares}. Then the natural functor
\[\mathbf{Qcoh}(X) \to \mathbf{Qcoh}(U) \times_{\mathbf{Qcoh}(U_Y)} \mathbf{Qcoh}(Y)\]
is an equivalence of categories. Consequently, the data of a vector bundle on $X$ is equivalent to the data of a vector bundle on $Y$, a vector bundle on $U$ along with a specified isomorphism between their restrictions to $U_Y$. More precisely, the natural functor
\[\mathbf{Vect}(X) \to \mathbf{Vect}(U) \times_{\mathbf{Vect}(U_Y)} \mathbf{Vect}(Y)\]
is an equivalence of categories.
\end{Theorem}

\begin{remark} \label{MVfaithful} To check that a vector bundle on $X$ is faithful (or $n$-twisted) one may pullback to each of the residual gerbes of $X$. Indeed, by definition, faithfulness (and $n$-twistedness) is a condition on the associated representation types. Thus, these can be checked on an open substack $U \subset X$ and its closed complement $Z=X \backslash U$. So even if $f$ (as in Definition \ref{DMVsquares}) is not representable, over the complement $Z=X \backslash U$, $f$ is an isomorphism by definition. Thus, if a sheaf on $X$ is faithful (or twisted) after restricting to $U$ and $Y$, it is so on $Z$ as well, and hence on $X$. We will combine this observation with Theorem \ref{MVglueing} to construct sheaves on stacks with desirable representation types. \end{remark}

A novelty of dealing with algebraic spaces and stacks is that points (and closed subsets) need not admit affine (or schematic) open neighborhoods, however they often admit affine \emph{flat} neighborhoods. Thus, in combination with the result above, such neighborhoods will give us a new way to produce vector bundles. Next, we will show that we can often \emph{construct} flat Mayer-Vietoris squares with favorable properties. The key technical tool is the theory of Tannakian duality first discussed in \cite{LurieTannaka} and developed further in \cite{hall2014coherent} and \cite{BhattTannaka}.

\begin{Theorem} \label{flatnhbdsexist} Let $X$ be a noetherian algebraic space. Suppose $Z=\Spec A_0 \to X$ is a closed immersion, then there is a flat neighborhood of $Z$ given by a noetherian affine scheme $Y$ and a map
\[f: Y=\Spec A \to X\]
with the property that $\dim(Y) \leq \dim(X)$ and $\dim(Y \backslash (Z \times_X Y)) \leq \dim(X)-1$.
\end{Theorem}

\begin{proof} Consider the quasi-coherent ideal sheaf $\cI \subset \cO_X$ corresponding to $Z$ and let $Z^{[i]}$ denote the closed subscheme corresponding to $\cI^{i+1} \subset \cO_X$. There is a sequence of closed immersions $Z^{[j-1]} \to Z^{[j]}$ for every $j \geq 1$. Since these are finite surjective morphisms and $Z$ is affine, \cite[Tag 07VP]{stacks} implies every $Z^{[i]}$ is an affine scheme. 

Since the global sections functor is exact on affine schemes we obtain for each $i \geq j$ a surjective ring map
\[\phi_{ij}: \Gamma(Z^{[i]}, \cO_X/\cI^{i+1}) = A_i \to  \Gamma(Z^{[i]}, \cO_X/\cI^{j+1})= A_{j}\]
If $I_i$ denotes the kernel of $A_i \to A_0$, then the kernel of $\phi_{ij}$ is $I_i^{j+1}$. Moreover, $I_1/I_1^2=I_1$ is a finitely generated $A_0$-module since it is an ideal of the noetherian ring $A_1$. It follows from \cite[EGA0 7.2.7]{grothendieck1960elements} that the limit $A=\lim A_i$ is an adic noetherian ring with respect to the kernel of the natural map $A \to A_0$, call it $I$.

By \cite[Corollary 1.5]{hall2014coherent} or \cite[Theorem 1.1]{BhattTannaka}, we obtain a map $f: Y= \Spec A \to X$, we will show that $Y$ is a flat neighborhood of $Z=\Spec A_0 \subset X$. First we check that $Y \times_X Z \to Z$ is an isomorphism. Note that this is certainly true at each finite stage i.e. 
\[Z^{[i]} \times_X Z \to Z\] 
is an isomorphism for every $i \geq 0$ (by construction of $Z^{[i]}$). Next, observe that $Y \times_X Z \to Y$ is a closed immersion and hence corresponds to a  $A$-algebra in $\Coh(\Spec A)$, call it $\mathscr{B}$. However, since $\Coh(\Spec A) \to \lim \Coh(\Spec A/I^i)$ is an equivalence of categories it suffices to determine the coherent algebra structure of $\mathscr{B}$ modulo $I^{i+1}$ for every $i \geq 0$. However, we have just seen that pulling back this coherent algebra to $\Spec A/I^{i+1}=Z^{[i]}$ yields the sheaf of algebras corresponding to the affine morphism $Z \to Z^{[i]}$. It follows that $Y \times_X Z \to Z$ is an isomorphism.

For each positive $i$, a similar argument shows that the projection map $Y \times_X Z^{[i]} \to Z^{[i]}$ is an isomorphism since $Z^{[j]} \times_X Z^{[i]} \cong Z^{[i]}$ over $X$ for all $j \geq i$. Since both $X$ and $Y$ are noetherian, a variant of the local criterion for flatness \cite[Tag 0523]{stacks} implies that $Y \to X$ is flat over the points of $Y \times_X Z=V(I) \subset \Spec A$. However, since $A$ is $I$-adic all maximal ideals must contain $I$ (see \cite[Theorem 8.2(i)]{Matsumura}) i.e. all closed points of $\Spec A$ are in $Z$. Thus, $f: \Spec A \to X$ is a flat morphism.

To see why $\text{dim}(Y) \leq \text{dim}(X)$ observe that at any closed point $x$ of $\Spec A$ the fiber $f^{-1}(f(x))$ is a singleton. Apply \cite[Theorem 15.1 (i)]{Matsumura} to conclude that the height of $x$ is no larger than the height of $f(x)$. Since $Z \times_X Y=Z \subset Y$ contains all the closed points of $Y$, it follows that $\dim(Y\backslash( Z \times_X Y)) \leq \dim(X)-1$. \end{proof}

\subsection{Rigidification of tame Artin stacks}

The purpose of this section is to prove that a normal tame Artin stack $\ms X$ which is finite-type and separated over a field $k$ can be rigidified. That is, there exists a morphism $\ms X \to \ms X^{rig}$ which gives $\ms X$ the structure of a gerbe over a stack $\ms X^{rig}$ with generically trivial stabilizers. In \cite{olssonboundedness} Olsson has shown this when $\ms X$ is a normal Deligne-Mumford stack over a locally Noetherian base. This implies our result when $k$ is of characteristic $0$, so it remains to show it for positive characteristics.

We work over the fppf site of a scheme $S$ and let $G$ and $P$ denote sheaves of groups on $S$. Following \cite[Chapter VII, pg. 134-144]{grothendieck2006groupes}, we will describe the category of group extensions of $P$ by $G$ using the language of bitorsors. By an extension of $P$ by $G$, we mean exact sequences of sheaves of groups
\[1 \to G \to E \to P \to 1\]
Here exactness means that for any $T \in S$ the sequence of groups is exact
\[1 \to G(T) \to E(T) \to P(T)\]
and that $E \to P$ is a surjective morphism of sheaves on $S$. If the above sequence of groups is exact, then the group structure on $E$ makes $E \to P$ into a $(G \times_S P, G \times_S P)$-bitorsor. Moreover, if $t_1,t_2 \in P(T)$ are two points, then multiplication in $E$ defines equivariant maps $\phi'_{t_1,t_2}: E_{t_1} \times E_{t_2} \to E_{t_1t_2}$ which are functorial in $T$. In the language of bitorsors: there is a isomorphism of $(G_T, G_T)$-bitorsors $\phi_{t_1t_2}: E_{t_1} \wedge E_{t_2} \to E_{t_1t_2}$ whose formation is functorial in $T$. Moreover, by the associativity of the multiplicative structure on $E$, these isomorphisms must satisfy the \emph{associativity relations}. That is, for all triples $t_1,t_2,t_3 \in P(T)$ and all $S$-schemes $T$, the following diagram must commute
\[\begin{tikzcd}
E_{t_1}\wedge E_{t_2} \wedge E_{t_3} \rar{\phi_{t_1,t_2} \wedge\text{ id }} \dar{\text{ id } \wedge \phi_{t_2,t_3}} & E_{t_1t_2} \wedge E_{t_3} \dar{\phi_{t_1t_2,t_3}} \\
E_{t_1} \wedge E_{t_2t_3} \rar{\phi_{t_1,t_2t_3}} & E_{t_1t_2t_3} \\
\end{tikzcd}\]
Thus, from the extension, we obtain the data $(G,P, E, \{\phi\})$ where $G$ and $P$ are sheaves of groups, $E$ is a $(G \times P, G \times P)$-bitorsor on $P$ and $\{\phi\}$ is a collection of isomorphisms as above. In fact, this data uniquely determines an extension of $P$ by $G$. 

Translating this into the setting where $G=\Delta$ is a diagonalizable group and $P=H$ is a finite constant group scheme: an extension of $H$ by $\Delta$ is equivalent to the data of $(\Delta, H, E, \{\phi\})$. However, we may simplify further because $\Delta$ is abelian. 

\begin{proposition} \label{bitorsors} The stack of $(\Delta,\Delta)$-bitorsors over $S_{fppf}$ is equivalent to $\mathrm{B}\Delta \times_S \Aut(\Delta)$. 
\end{proposition}

The proof is standard so we omit it. To summarize everything, we have:

\begin{proposition} \label{extensionreductive} The category of extensions of a constant group by a diagonalizable group is equivalent to the category of tuples $(H,\Delta, (E,\psi), \{\phi\})$ where 

\begin{enumerate}
\item $H$ is a constant group scheme on $S$.
\item $\Delta$ is a diagonalizable group scheme on $S$.
\item $E \to H$ is a left $\Delta \times_S H$-torsor, and $\psi \in \mathrm{Aut}(\Delta) (H)$, an automorphism of $\Delta \times_S H$.
\item For every pair $h_i,h_j \in H(S)$, an isomorphism $\phi_{h_ih_j}: E_{h_i} \wedge E_{h_j} \to E_{h_ih_j}$ of bitorsors satisfying the associativity relations.
\end{enumerate} \end{proposition}

\begin{lemma} \label{connected} Let $T$ be a noetherian scheme over a field $k$ of positive characteristic $p$ and suppose $G \to T$ is a finite flat linearly reductive group scheme of length $p^nl$ where $(p,l)=1$. Then the subfunctor $\Delta \subset G$ whose $B$-points are defined as 
\[\Delta(B)=\{g \in G(B)|\exists m\geq 0\text{ such that } g^{p^m}=1\}\]
is representable by a finite flat $T$-scheme. Moreover, $\Delta$ is a locally diagonalizable group scheme which is normal in $G$ and fits into an exact sequence of finite flat group schemes over $T$
\[0 \to \Delta \to G \to G/\Delta \to 0\]
where $G/\Delta$ is tame and \'etale.
\end{lemma}

\begin{proof} Since $\Delta \subset G$ is a subsheaf and the statements are all local, to show the lemma we may check it locally on $T$. Thus, we assume that $k$ is algebraically closed. 

First consider the case when $T=k$ and observe that $\Delta$ coincides with the functor corresponding to the connected component $G_0$ of the identity of $G$. Indeed, the quotient $G/G_0$ is \'etale (hence constant) and since it is also linearly reductive the order of $G/G_0$ is prime to $p$. It follows that $\Delta \subset G_0$. Note that the connected component $G_0$ has order a power of $p$ by \cite[3.7(II)]{tate1997finite} and because a finite group scheme over a field is killed by its order (see \cite[3.7(III)]{tate1997finite}), we have $G_0 \subset \Delta$. Thus, the lemma holds when $T=\Spec k$ by \cite[Lemma 2.11(a)]{Tame}.
 
The general case follows from the structure theory of finite flat linearly reductive group schemes. Indeed, there is a collection of jointly surjective morphisms $T_i \to T$ and an isomorphism of group schemes $G|_{T_i} \cong (H \ltimes \Delta') \times_k T_i$ where $H$ is constant and tame and $\Delta'$ is diagonalizable (see \cite[Theorem 2.19]{Tame}). In other words, $G|_{T_i}$ is the pullback of a well-split group scheme on $\Spec k$ and the subfunctor $\Delta'|_{T_i} \subset G|_{T_i}$ is equal to the (subfunctor associated to the) pullback of the connected component: $(H \ltimes \Delta')_0 \subset H \ltimes \Delta'$. The result now follows from the case $T=\Spec k$. \end{proof}

\begin{lemma} \label{extendingsubgroup} Fix a field $k$ of positive characteristic $p$ and let $R$ be a strictly local Henselian $k$-algebra which is noetherian and normal. Let $V \subset \Spec R$ denote a dense open subset of $\Spec R$. Suppose that $G \to \Spec R$ is a finite flat linearly reductive group scheme, and $G'_V \subset G|_V$ is a closed subgroup scheme which is also finite flat over $V$. Then there is a finite flat subgroup $G' \subset G$ over $\Spec R$ with equality $G'|_V=G'_V$ as subgroup schemes of $G|_V$. 
\end{lemma}

\begin{proof} By Lemma \ref{connected} applied to $G \to \Spec R$ the subfunctor $\Delta \subset G$ consisting of elements of order a power of $p$ induces an extension of finite flat linearly reductive group schemes
\[0 \to \Delta \to G \to H \to 0\]
where $\Delta$ is locally diagonalizable and $G/\Delta$ is tame and \'etale. Since $\Spec R$ is strictly Henselian, it follows that $\Delta$ is diagonalizable and $H$ is constant. Upon restriction to $V$ and applying Lemma \ref{functorial} to the inclusion $G'_V \subset G|_V$ we obtain an extension
\[0 \to \Delta' \to G'_V \to H' \to 0\]
where $\Delta'$ is diagonalizable and $H'$ is constant. Thus, $\Delta'$ and $H'$ both extend to group schemes over $\Spec R$. We will show that the extension also extends by showing that the corresponding $(\Delta', \Delta')$-bitorsor over $H'$ extends. However by Proposition \ref{extensionreductive} we may forget the bitorsor structure if we remember the underlying left-torsor along with the automorphism $\Delta' \to \Delta'$ specifying the right-action. Since $\text{Aut}(\Delta')$ is constant, any automorphism \emph{uniquely} extends over $\Spec R$. Thus, to extend the bitorsor $G'_V \to H'$ and the bi-equivariant inclusion $G'_V \subset G|_V$, it suffices to extend the underlying left $\Delta'$-torsor and the left-equivariant inclusion.

The inclusion $G'_V \subset G|_V$ is equivariant for the map $\Delta' \to \Delta$. In other words, if we extend the structure group of the left $\Delta'$-torsor $G'_V \to H'$ along the map $\Delta' \to \Delta$, we obtain the left $\Delta$-torsor $G|_V \times_H H' \to H'$. However $G|_V \to H|_V$ is the restriction of the map $G \to H$ over $\Spec R$, i.e. the left $\Delta$-torsor $G|_V \times_H H' \to H'$ extends to $\Spec R$. This implies the left $\Delta'$-torsor $G'_V \to H'$ extends to $\Spec R$ as well. Indeed, there is an exact sequence of diagonalizable groups (over $\Spec R$)
\[0 \to \Delta' \to \Delta \to \Delta'' \to 0\]
and consider the long exact sequences induced by the functors $\H^0(\Spec R, -)$ and $\H^0(V,-)$:
\[\begin{tikzcd}
\H^0(\Spec R, \Delta'') \arrow[r] \arrow{d}{\cong} & \H^1(\Spec R, \Delta') \arrow[d,"\phi'", hook] \arrow[r,"g"] & \H^1(\Spec R, \Delta) \arrow[d,"\phi", hook] \arrow[r] & \H^1(\Spec R, \Delta'') \arrow[d, hook]  \\
\H^0(V, \Delta'') \arrow[r] & \H^1(V, \Delta')  \arrow[r,"f"] & \H^1(V, \Delta) \arrow[r] & \H^1(V, \Delta'') \\
\end{tikzcd}\]
The left-most vertical arrow is injective because $V \subset \Spec R$ is schematically dense and $\Delta''$ is separated. Surjectivity of the left-most arrow and injectivity of the remaining arrows follows from the normality of $\Spec R$. Indeed, the closure of a $V$-point in a finite $\Spec R$-scheme is a $\Spec R$-point because the closure is finite and birational over $\Spec R$ (see \cite[Tag 0AB1]{stacks}). Now apply this either to the finite scheme $\Delta''$ to deduce the surjectivity of the left-most arrow or a left-torsor under a finite group scheme to deduce the injectivity of the remaining three arrows.

A quick diagram chase shows that if $x \in \H^1(V, \Delta')$ has the property that $f(x)=\phi(y)$ for $y \in \H^1(\Spec R, \Delta)$, then there is a $y' \in \H^1(\Spec R, \Delta')$ with $\phi'(y')=x$ and $g(y')=y$. Thus, both the underlying left $\Delta'$-torsor $G'_V \to H'$ and its equivariant inclusion into $G|_V$ extends over $\Spec R$ and hence the $(\Delta',\Delta')$-bitorsor does as well.

It remains to show that for every $h_i,h_j \in H'(\Spec R)$, there are isomorphisms $\phi_{h_ih_j}: G'_{h_i} \wedge G'_{h_j} \to G'_{h_ih_j}$ of $(\Delta',\Delta')$-bitorsors satisfying the associativity relations. In this case, the data of an isomorphism of bitorsors is the same as the data of an isomorphism of the underlying left torsors. Indeed, a left-equivariant map is right-equivariant exactly when the two sections of $\text{Aut}(\Delta')$ determining the right action coincide. However, this automorphism scheme is constant and since $\phi_{h_ih_j}$ exists over $V$, any such isomorphism of left $\Delta'$-torsors is automatically an isomorphism of bitorsors. In other words, we desire sections of $\text{Isom}_{\Delta'}(G'_{h_i} \wedge G'_{h_j}, G_{h_ih_j}')$  defined over $\Spec R$. But this sheaf is itself a left $\Delta'$-torsor and is therefore finite over $\Spec R$. Thus, since we had such maps $\phi_{h_ih_j}$ over the open subset $V \subset \Spec R$, the closure of the corresponding $V$-points in the $\text{Isom}$-scheme, must be $\Spec R$-points. That this collection of isomorphisms satisfy the associativity relations follows because they are satisfied over $V$. The result follows. \end{proof}

\begin{lemma} \label{functorial} Let $T$ be scheme finite-type over a field $k$ of positive characteristic $p$ and suppose that $G \to T$ is a finite flat linearly reductive group scheme with length $n$. Suppose that the subgroup schemes $\Delta$ and $G/\Delta$ defined in Lemma \ref{connected} are diagonalizable and constant respectively. Then if $G' \subset G$ is any finite flat subgroup scheme over $T$, it appears in an exact sequence
\[0 \to \Delta' \to G' \to G'/\Delta' \to 0\]
where $\Delta'$ is diagonalizable and $G'/\Delta'$ is constant and tame. 
\end{lemma} 

\begin{proof} Note that $G'$ is linearly reductive since it is a subgroup of $G$. Next, define $\Delta'=\Delta \times_G G' \subset G'$. Note that $\Delta'$ coincides with the subgroup of elements of $p$-power order in $G'$, thus it is locally diagonalizable, normal in $G'$, and the quotient $G'/\Delta'$ is \'etale and tame by Lemma \ref{connected}. However, $\Delta' \subset \Delta$ and since the latter is diagonalizable, $\Delta'$ is as well. The induced map $G'/\Delta' \to G/\Delta$ must be injective. Indeed, if $K$ denotes the kernel, then one can (tediously) check that there is an injective morphism of group schemes $K \to \Delta/\Delta'$ over $T$. If $G$ was abelian, this would follow directly from the snake lemma but in our case the fact that $\Delta$ is abelian suffices. Now because $K$ is finite \'etale with geometrically connected fibers over $T$ (since it is a subscheme of the locally diagonalizable scheme $\Delta/\Delta'$ whose topological space is a single point), it must be the trivial group scheme. It follows that $G'/\Delta' \subset G/\Delta$ and because the latter is constant and tame, the former is as well. \end{proof}

\begin{Theorem} \label{tamerigidification} Let $\ms X$ be a normal tame Artin stack which is separated and finite-type over a field $k$ of positive characteristic $p$, then the morphism to the coarse moduli space of $\ms X$ factorizes 
\[\ms X \to \ms X^{rig} \to X\]
where $\ms X$ is a gerbe over $\ms X^{rig}$ and the latter has generically trivial stabilizers. 
\end{Theorem}

\begin{proof} Our goal will be to produce a finite flat subgroup stack $G \subset I_{\ms X}$ containing the generic fiber of $I_{\ms X} \to \ms X$. This suffices because we may apply \cite[Theorem A.1]{Tame} to produce the desired rigidification. 

Since $\ms X$ is reduced there is a open substack $\ms U \subset \ms X$ where $I_{\ms X}|_{\ms U} \to \ms U$ is flat. Let $G \subset I_{\ms X} \to \ms X$ denote the scheme-theoretic image of the map $I_{\ms X}|_{\ms U} \to I_{\ms X}$, note that it is closed by construction. We claim that this is a flat subgroup stack of $I_{\ms X}$. The scheme-theoretic image, $j: G \to I_{\ms X}$, contains the identity section since it contains the schematically dense open subset $\ms U$ of $\ms X \subset I_{\ms X}$. Moreover, it is closed under multiplication $m: G \times_{\ms X} G \to I_{\ms X}$ and inversion $i: G \to I_{\ms X}$ if and only if the following cartesian diagrams

\begin{minipage}{.5\textwidth}
\[\begin{tikzcd}
K \rar{f} \dar & G \times_{\ms X} G \dar{m} \\
G \rar{j} & I_{\ms X} \\
\end{tikzcd}\]
\end{minipage}
\begin{minipage}{.5\textwidth}
\[\begin{tikzcd}
K' \rar{f'} \dar & G \dar{i} \\
G \rar{j} & I_{\ms X} \\
\end{tikzcd}\]
\end{minipage}

\noindent have the property that $f$ and $f'$ are isomorphisms. Note that $f'$ is an isomorphism because $K'$ contains $I_{\ms X}|_{\ms U}$ which is schematically dense in $G$. On the locus $\ms V \subset \ms X$ over which $G$ is flat (so $\ms U \subset \ms V$), it is also true that $I_{\ms X}|_{\ms U} \times_{\ms X} I_{\ms X}|_{\ms U}$ is schematically dense in $G|_{\ms V} \times_{\ms X} G|_{\ms V}$ (see \cite[Tag 01RE]{stacks}). Thus, because $K$ contains this schematically dense substack over $\ms V$, it follows that $f$ is an isomorphism over $\ms V$.

It follows that $G$ is a subgroup stack of $I_{\ms X}$ over $\ms X$ if it is flat over $\ms X$, and this is what we show next. Since the formation of a scheme-theoretic image commutes with flat base change on the target (see \cite[Tag 0CMK]{stacks}), we may check $G$ is flat over $\ms X$ fpqc-locally on $\ms X$. Thus, we may pass to the strict Henselization $\Spec R \to X$ of a closed point of the coarse moduli space $X$ of $\ms X$. By \cite[Theorem 3.2(d)]{Tame} we have $\ms X=[Y/F]$ where $Y$ is a finite $\Spec R$-scheme and $F$ is a locally well-split group scheme over $\Spec R$. Let $W \to \ms X$ be a flat cover by a normal scheme where $W=\Spec R'$ is the spectrum of a strictly local ring and set $U=\ms U \times_{\ms X} W$. We now have the following representable morphisms:
\[W \to [Y/F] \to \text{BF}\]
and pulling back the inertia stacks yields the following closed immersions over $W$:
\[G \times_{\ms X} W \subset I_{\ms X} \times_{\ms X} W \subset I_{\text{BF}} \times_{\text{BF}} W=I\]
where $G \times_{\ms X} W$ is the scheme-theoretic closure of $G \times_{\ms X} U=G_U$ in $I_{\ms X} \times_{\ms X} W$ and $I$ is an (inner) form of $F \times_{\Spec R} W$. Hence $I$ is a finite flat linearly reductive group scheme over $W$.

Since $G_U$ is a finite flat subgroup scheme of $I|_U$, Lemma \ref{extendingsubgroup} implies there is a finite flat closed subgroup scheme $G' \subset I$ over $W$ which extends the inclusion $G_U \subset I|_U$. Note that the open immersion $G_U=G'|_U \subset G'$ is schematically dense and therefore $I_{\ms X} \times_{\ms X} W$ must contain $G'$ as a closed subgroup scheme. It follows that the closure of $G_U$ in $I_{\ms X} \times_{\ms X} W$ is $G'$. In particular, the closure of $G_U$ in $I_{\ms X} \times_{\ms X} W$ is finite flat over $W$, as desired. \end{proof}

\section{The Resolution Property of Algebraic Space Surfaces and Tame Artin Stacks}

\subsection{Separated Algebraic Space Surfaces have the Resolution Property}

In this section, we prove the following

\begin{Theorem} \label{ResSurface} Let $X$ be a $2$-dimensional algebraic space which is finite-type and separated over a noetherian ring $R$. If either
\begin{enumerate} 
\item $X$ is catenary and Jacobson, or
\item  the morphism $X \to \Spec R$ is proper,
\end{enumerate}
then $X$ has the resolution property.

 \end{Theorem}

The novelty of Theorem \ref{ResSurface} is that it addresses the first class of algebraic spaces which admit non-schematic examples. It should be noted that there exist examples of non-schematic algebraic spaces known to have the resolution property (see, for example, the discussion following \cite[Theorem 2.1]{TotaroRes}). However, the scope of these examples is fairly narrow: they all admit finite flat covers by schemes. Even for normal algebraic space surfaces, such covers are difficult to construct.

\begin{remark} \label{spacetostack} Unfortunately the methods used to prove this result cannot be used to remove the normality hypothesis in Theorem \ref{T1}. For instance, here we repeatedly use the fact that algebraic spaces are generically schemes, whereas non-normal tame Artin stacks cannot always be rigidified, so I am not aware of a way to apply these methods to the general stacky case. \end{remark}

The schematic case was treated by Schr\"{o}er-Vezzosi \cite{NormalSurfaceRes} and Gross' \cite{SurfaceRes}. There are three serious reasons why their methods do not extend to the algebraic space setting. 
\begin{enumerate}
\item The construction of an "almost-ample" family of rank $1$ sheaves (see \cite[Proposition 1.8]{SurfaceRes} and \cite[Proposition 2.2]{NormalSurfaceRes}) does not work for algebraic spaces. Indeed, for this one requires that every point of $X$ admit an affine open neighborhood. This reduces the problem to resolving sheaves which are locally free away from a finite set of points. We make a similar reduction by using formal-local methods to construct higher rank resolutions (see Lemma \ref{MVresolution}). 
\item Moreover, the "almost-ample" family of rank $1$ sheaves have positivity properties that our higher rank sheaves do not. To remedy this we use Lemma \ref{tensortracetrick} and a refinement of Gross' construction of a $1$-ample family of vector bundles.

\item Gross believed that his local analysis (see \cite[Proposition 2.7]{SurfaceRes}) required the resolution property \emph{Zariski}-locally. However, using flat descent as in \cite{hall2016mayer}, we argue that it suffices to work over flat neighborhoods. In this way, we establish a similar cohomological obstruction to resolving sheaves. 

\end{enumerate}

\begin{definition} \label{FK} (\cite[2.2]{SurfaceRes}) Let $X$ be a locally noetherian algebraic space, we say that a coherent sheaf $F$ is $F_k$ if $F|_{\Spec \mathcal{O}^h_{X,x}}$ is free for all $x \in X$ with $\text{dim}\mathcal{O}^{h}_{X,x} \leq k$ and for all other $x$ we have $\text{pd}_{\mathcal{O}^{h}_{X,x}}(F|_{\Spec \mathcal{O}^h_{X,x}}) \leq \text{dim} \mathcal{O}^h_{X,x} - k$.

\end{definition}

\noindent First, we describe a cohomological obstruction to the existence of a resolution by $F_1$ sheaves. We use flat neighborhoods and formal glueing to extend \cite[Proposition 2.7]{SurfaceRes} to the setting of algebraic spaces.

\begin{proposition} \label{flat-glueing} Let $X$ be a $2$-dimensional, noetherian, algebraic space and $F$ a coherent sheaf which is locally free of constant rank away from a closed set $Z \subset X$ of codimension $2$ which contains no associated points of $X$. There exists
\begin{enumerate}
\item a coherent sheaf $L$ on $X$ with $L|_{X-Z} \cong \det(F|_{X-Z})^{\vee}$, and
\item an obstruction $o \in \H^2(X, \underline{\mathrm{Hom}}(F,L))$, whose vanishing is sufficient to guarantee the existence of an exact sequence of coherent sheaves on $X$:
\[0 \to L \to N \to F \to 0\]
where $N$ is $F_1$.
\end{enumerate}

\end{proposition}

\begin{proof} Let $Y=\bigsqcup_{z \in Z} \Spec R_z$ be a disjoint union formed by the Henselian local rings $(R_z,m_z)$ of each $z \in Z$ (see \cite[Tag 0BGW]{stacks}). Observe that the natural map $Y \to X$ is a flat neighborhood of $Z$. Therefore if we set $U=X\backslash Z$, we obtain a flat Mayer-Vietoris square.

\[\begin{tikzcd}
U_Y=\bigsqcup_{z \in Z} (\Spec R_z \backslash \{m_z\}) \rar \dar & Y=\bigsqcup_{z \in Z} \Spec R_z \dar{f} \\
U \rar{i} & X \\
\end{tikzcd}\]

Restrict $F$ to $Y$ and on each noetherian local ring $R_z$ apply \cite[Lemma 2.6]{SurfaceRes} to $F|_{\Spec R_z}$. This yields an exact sequence, or an extension class, call it $\sigma_z$:
\[0 \to L_z \to N_z \to F|_{\Spec R_z} \to 0\]
of $R_z$ modules where $N_z$ and $L_z$ are locally free away from $m_z$, $N_z$ has projective dimension $ \leq 1$ and isomorphisms
\[L_z|_{(\Spec R_z \backslash \{m_z\})} \cong \det (F|_{(\Spec R_z \backslash \{m_z\})})^{\vee} \cong \det (F|_U)^{\vee}|_{(\Spec R_z \backslash \{m_z\})}\]
These isomorphisms and the formal glueing Theorem \ref{MVglueing} imply that there exists a coherent sheaf $L \in \Qcoh(X)$ which restricts to $\det(F|_U)^{\vee}$ over $U$ and $\bigsqcup L_z$ over $Y$.

As such, we may view the sequences above as extension classes $\sigma_z \in \text{Ext}^1(F|_{\Spec R_z}, L_z)$ and we claim these glue to yield a global section of $E=\underline{\text{Ext}}^1(F,L)$. Indeed, a global section of $E$ corresponds to a morphism of $\mathcal{O}_X$-modules: $\mathcal{O}_X \to E$. Moreover, we have $E|_U=0$ (because $F$ is locally free on $U$) and the map $(\sigma_z)_{z \in Z}: \mathcal{O}_Y \to E|_Y=\text{Ext}^1(F|_Y,L|_Y)$. We may view this data as a morphism of triples
\[(0,(\sigma_z)_{z \in Z}): (\mathcal{O}_U, \mathcal{O}_Y, \text{id}) \to (0, \underline{\text{Ext}}^1(F|_Y, L|_Y), 0)\]
i.e. a morphism in the category $\Qcoh(U) \times_{\Qcoh(U_Y)} \Qcoh(Y)$. Thus, by Theorem \ref{MVglueing}, there is a section $\sigma \in \H^0(X, \underline{\text{Ext}}^1(F,L))$ that extends $(\sigma_z)_{z \in Z}$. Finally, by considering the low degree terms of the local-global Ext spectral sequence
\[\text{H}^p(X,\underline{\text{Ext}}^q(F,L)) \Rightarrow \text{Ext}^{p+q}(F,L)\]
we obtain the exact sequence
\[\text{Ext}^1(F,L) \to \text{H}^0(X, \underline{\text{Ext}}^1(F,L)) \to \text{H}^2(X, \underline{\text{Hom}}(F,L))\]
Therefore, the image of $\sigma$ in $\text{H}^2(X, \underline{\text{Hom}}(F,L))$ is the desired obstruction. Indeed, the vanishing of $\sigma$ is equivalent to the existence of a global extension
\[0 \to L \to N \to F \to 0\]
extending the local extension classes $(\sigma_z)_{z \in Z}$. \end{proof}

Next, we describe a cohomological obstruction to resolving $F_1$ sheaves by vector bundles. The following extends \cite[Proposition 2.4]{SurfaceRes} to the setting of algebraic spaces.

\begin{proposition} \label{Gross2} Let $X$ be a dimension $2$, noetherian, algebraic space and $F$ a coherent sheaf which is $F_1$. Let $Z$ denote the (finite) locus where $F$ is not locally free.  Assume that $F$ has constant rank at all generic points of $X$. Fix a vector bundle $V$ on $X$, then for any $m \gg 0$, there is an obstruction $\sigma \in \H^2(X, \underline{\mathrm{Hom}}(F, V^{\oplus m}))$ whose vanishing guarantees the existence of a locally free resolution
\[0 \to V^{\oplus m} \to E \to F \to 0\]
where $E$ is locally free.
\end{proposition}

\begin{proof} The proof is very similar to the previous proposition, so we just give a sketch here. Let $Y = \bigsqcup_{z \in Z} \Spec R_z \to X$ be a flat neighborhood of $Z$ which is the disjoint union of the Henselian local rings of the $z \in Z$. Since $F$ is $F_1$ and $Y$ is semi-local, there exists a large integer $m$ and an exact sequence of sheaves
\[0 \to V^{\oplus m}|_Y \to \mathcal{O}_Y^{\oplus n} \to F|_Y \to 0\]
Thus we obtain a section $\sigma_Y \in \H^0(Y, \underline{\text{Ext}}^1(F|_Y, V^{\oplus m}|_Y))$. Just as in the proof of the previous proposition we can use Theorem \ref{MVglueing} to infer the existence of a global section $\sigma$ of $\underline{\text{Ext}}^1(F, V^{\oplus m})$ which extends $\sigma_Y$. Using the low degree terms of the local-global Ext spectral sequence, the obstruction to lifting $\sigma \in \H^0(X, \underline{\text{Ext}}^1(F,V^{\oplus m}))$ to an element in $\text{Ext}^1(F, V^{\oplus m})$ lies in $\H^2(X, \underline{\text{Hom}}(F, V^{\oplus m}))$. The result follows. \end{proof}

Next we show that any coherent sheaf on a $2$-dimensional algebraic space (with affine diagonal) can be resolved by a sheaf $F$ which is locally free away from a finite set of points. First, we record the following

\begin{lemma} \label{cohomologicaldimension1} Let $X$ be an algebraic space with affine diagonal and an open subset $j: \Spec A \to X$ with the property that $(X \backslash \Spec A)_{\text{red}}$ is affine. Then $\H^i(X,F)=0$ for every $i \geq 2$ and every quasi-coherent sheaf $F$ of finite-type. 
\end{lemma} 

\begin{proof} Fix a quasi-coherent sheaf $F$ on $X$. From the natural map $f: F \to j_*j^*F$, we obtain the following exact sequences
\[0 \to \text{Ker}(f) \to F \to \text{Im}(f) \to 0\]
\[0 \to \text{Im}(f) \to j_*j^*F \to \text{Coker}(f) \to 0\]
By hypothesis, the sheaves $\text{Coker}(f)$ and $\text{Ker}(f)$ are supported on an affine closed subscheme. Moreover, because $X$ has affine diagonal we have
\[\H^i(X,j_*j^*F)=\H^i(\Spec A, j^*F)=\H^i(X, \text{Ker}(f))=\H^i(X, \text{Coker}(f))=0\]
for every $i \geq 1$. It follows that
\[\H^i(X, F)=\H^i(X,\text{Im}(f))\]
for $i \geq 1$ and the latter vanishes for every $i \geq 2$. \end{proof}

\begin{lemma} \label{freeaffine} Let $X$ be a noetherian algebraic space with affine diagonal and $j: \Spec A \to X$ an open immersion. Every coherent sheaf $M$ on $X$ admits a surjection $F \to M \to 0$ where $F$ is a coherent sheaf which is free over $\Spec A$. 
\end{lemma}

\begin{proof} Since $\Spec A$ is affine, there exists a surjection $V \to j^*M \to 0$ where $V$ is a free module. Moreover, because $j$ is an affine map, we obtain a surjection $j_*V \to j_*j^*M \to 0$. Let $K$ denote the fiber product of this map along the natural morphism $M \to j_*j^*M$. This yields a surjection $K \to M \to 0$ where $K$ is quasi-coherent and $K|_{\Spec A} \cong V$ is free. However, because $M$ and $K|_{\Spec A}$ are coherent, and $K$ can be written as a colimit of its coherent subsheaves (see \cite[Proposition 15.4]{champs}), there is a coherent subsheaf $K' \subset K$ with $K'|_{\Spec A} \cong K|_{\Spec A}$ and where the induced map $K' \to M$ is surjective. \end{proof}

\begin{lemma} \label{MVresolution} Let $X$ be a noetherian $2$-dimensional algebraic space with affine diagonal. Suppose that $j: \Spec A \to X$ is an open subscheme and $(X \backslash \Spec A)_{\mathrm{red}}$ is affine. Then $X$ satisfies the resolution property.
\end{lemma}

\begin{proof} Fix a coherent sheaf $M$, by Lemma \ref{freeaffine} there is a coherent sheaf $F$ which is locally free over $\Spec A$ and a surjection $F \to M \to 0$. Thus, by replacing $F$ with $M$, we may assume $M$ is locally free over $\Spec A$. Since $(X \backslash \Spec A)_{\text{red}}=\Spec B_0 \to X$ is an affine closed subscheme, Theorem \ref{flatnhbdsexist} implies there is a flat Mayer-Vietoris square
\[\begin{tikzcd}
\Spec D \rar \dar & \Spec B \dar{f} \\
\Spec A \rar{i} & X \\
\end{tikzcd}\]
where $\dim(\Spec D) \leq 1$. Now, on $\Spec B$ there exists an extension, call it $\sigma_B$: 
\[0 \to N' \to \mathcal{O}_{\Spec B}^{\oplus n} \to M|_{\Spec B} \to 0\]
and because $M|_{\Spec D}$ is locally free, it follows that $N'|_{\Spec D}$ is as well. Taking determinants over $\Spec D$, we obtain an isomorphism $\det(M|_{\Spec D})^{\vee} \cong \det(N'|_{\Spec D})$. Moreover, by Serre's splitting theorem (see \cite[Theorem 1]{SerreSplitting}) we have isomorphisms
\[N'|_{\Spec D} \cong \mathcal{O}_{\Spec D}^{\oplus m-1} \oplus \det(N'|_{\Spec D}) \cong \mathcal{O}_{\Spec A}|_{\Spec D}^{\oplus m-1} \oplus \det(M|_{\Spec A})^{\vee}|_{\Spec D}\]
and by Theorem \ref{MVglueing} this implies there is a coherent sheaf on $X$ extending $N'$ which is locally free on $\Spec A$, call this sheaf $N$. Following the proof of Proposition \ref{flat-glueing}, the exact sequence above yields a global section $\sigma \in \H^0(X, \underline{\text{Ext}}^1(M,N))$ which restricts to $\sigma_B$ over $\Spec B$. By Lemma \ref{cohomologicaldimension1} the local-global Ext-spectral sequence tells us $\sigma$ comes from a class in $\text{Ext}^1(M,N)$. In other words, there is an extension on $X$:
\[0 \to N \to F' \to M \to 0\]
where $F'$ is locally free. Indeed, over the open subset $\Spec A$ it lies in between two locally free sheaves and hence is locally free as well. Over the points of $X \backslash \Spec A$, we may check this after passing to the flat cover $\Spec B \to X$. Here the sequence above restricts to the extension $\sigma_B$ so $F'|_{\Spec B} \cong \mathcal{O}_B^{\oplus n}$. Thus, $X$ has the resolution property. \end{proof}

\begin{proposition} \label{freethick} Let $X$ be a noetherian $2$-dimensional algebraic space with affine diagonal. Then every coherent sheaf $M$ admits a surjection from a coherent sheaf $F \to M$ where $F$ is locally free over a thick open subset $U \subset X$ (i.e. the complement of $U$ is of codimension at least $2$).
\end{proposition}

\begin{proof} Let $\Spec A \subset X$ be a dense open subscheme. The (reduced) complement of $\Spec A$, call it $C$, admits a finite affine open covering $C=\bigcup_{i=1}^n \Spec B_i$ where each $\Spec B_i \subset C$ is dense. Since $C$ has the subspace topology (see \cite[Tag 04CD]{stacks}), there are open subsets $W_i \subset X$ with $|W_i \cap C|=|\Spec B_i|$. Thus, the open subsets $f_i: W_i \cup \Spec A \to X$ are thick open neighborhoods which cover $X$. Moreover, the reduced complement $((W_i \cup \Spec A) \backslash \Spec A)_{\text{red}}$ is supported on the affine closed subscheme $\Spec B_i \to W_i \cup \Spec A$. Thus, each $W_i \cup \Spec A$ satisfies the hypothesis of Lemma \ref{MVresolution}, so for each $i$ there exists a vector bundle $F'_i$ on $W_i \cup \Spec A$ and a surjection $F'_i \to M|_{W_i \cup \Spec A}$. 

Take the pushforward $(f_i)_*F'_i \to (f_i)_*(f_i)^*M$ and define $F''_i$ to be the pullback along the natural map $a: M \to (f_i)_*(f_i)^*M$. This yields a map from a quasi-coherent sheaf $F''_i \to M$ which is surjective over $W_i \cup \Spec A$ and where $F''_i$ is locally free on $W_i \cup \Spec A$. Thus $F''=\bigoplus_{i=1}^n F''_i \to M$ is surjective and $F''$ is locally free on the thick open neighborhood $U= \bigcap_{i=1}^n (W_i \cup \Spec A)$. Since $M$ and the $F''|_U$ are coherent, there is a coherent subsheaf $F \subset \bigoplus_{i=1}^n F''_i$ which surjects onto $M$ and which is locally free on $U$. \end{proof}

\begin{definition} (\cite[3.1]{SurfaceRes}) Consider a family of vector bundles $\{\mathcal{E}_n\}_{n \in \mathbf{N}}$ on an algebraic space $X$ of dimension $m$ and let $d \geq m-1$. Then $\{\mathcal{E}_n\}_{n \in \mathbf{N}}$ is said to be a (cohomologically) $d$-ample family of vector bundles if for every coherent sheaf $M$ there is a $n_0 \geq 0$ such that $\H^i(X, M \otimes \mathcal{E}_n)=0$ for all $n \geq n_0$ and $i \geq d+1$.
\end{definition}

In fact, every $2$-dimensional algebraic space proper over a noetherian ring admits a family of $1$-ample vector bundles of rank $2$. The proof is essentially the same as for schemes (see \cite[Theorem 4.5]{SurfaceRes}) but there are a few modifications required, we indicate them below. 

\begin{lemma} \label{finitedescent} Let $f: X' \to X$ be a finite morphism of noetherian algebraic spaces which is an isomorphism over a dense open subset $U \subset X$. If $X \backslash U$ consists of a finite set of points, then the pullback functor $f^*: \textbf{Vect}(X) \to \textbf{Vect}(X')$ is essentially surjective. 
\end{lemma}

\begin{proof} There is a flat neighborhood $Y=\bigsqcup_{z \in X \backslash U} \Spec R_z \to X$ formed by the spectra of the Henselian local rings $R_z$ and let $Y^{\circ}=Y \times_X U$. Pulling back the resulting flat Mayer-Vietoris square along the map $X' \to X$ yields another flat Mayer-Vietoris square. Set $Y'=Y \times_X X'$ and note that $Y^{\circ} \times_X X'=Y^{\circ}$ since $f$ is an isomorphism over $U$. Moreover, by Proposition \ref{MVglueing} to show the functor $f^*$ is essentially surjective is suffices to show the natural functor
\[\textbf{Vect}(U) \times_{\textbf{Vect}(Y^{\circ})} \textbf{Vect}(Y) \to \textbf{Vect}(U) \times_{\textbf{Vect}(Y^{\circ})} \textbf{Vect}(Y')\]
is essentially surjective. Since $Y$ is semi-local, so is $Y'$ by \cite[Tag 04GG (9)]{stacks} and therefore there are no nontrivial vector bundles on $Y'$. It follows that any object on the right hand side is in the essential image. \end{proof}

\begin{lemma} \label{1ample} (Gross) Fix a Noetherian ring $R$ and let $X$ be a $2$-dimensional algebraic space which is proper over $\Spec R$. Then $X$ admits a cohomologically $1$-ample family $\{\mathcal{E}_n\}$ of rank $2$ vector bundles so that $\det(\mathcal{E}_n)|_U=\mathcal{O}_X|_U$ where $\mathrm{codim}(X \backslash U) \geq 2$. 
\end{lemma} 

\begin{proof} The proof proceeds by taking a proper birational morphism from a projective scheme $f: Y \to X$ which is an isomorphism over a thick open subset $U \subset X$ and constructing a $1$-ample family of rank $2$ vector bundles $\{\mathcal{F}_n\}$ on $Y$. By ensuring that the $\mathcal{F}_n$ are trivial on a fixed thickening of the exceptional locus of $f$, one may conclude that $\mathcal{F}_n$ descends to a vector bundle $\mathcal{E}_n$ on $X$ by \cite[Lemma 4.1]{SurfaceRes} (see also \cite[Proposition 1.2]{NormalSurfaceRes}).

To obtain a thick quasiprojective open subset $U \subset X$ use \cite[Tag 0ADD]{stacks} and \cite[Theorem 1.5]{SurfaceRes}. Then, to construct a $U$-admissible blow-up $f: Y \to X$ where $Y$ is quasiprojective apply \cite[Tag 088R]{stacks}. The rest of the proof follows exactly as in \cite[Theorem 4.5]{SurfaceRes} except for a minor detail in \cite[Lemma 4.1]{SurfaceRes}. To show that $\mathcal{F}_n$ descends along the morphism $Y \to X$, the proof in \cite[Proposition 1.2]{NormalSurfaceRes} first considers the finite morphism in the Stein factorization $Y \to X' \to X$. The argument there uses the existence of affine open neighborhoods but can be modified easily. Indeed, the map $X' \to X$ is finite and an isomorphism away from a finite set of points. It follows that the pullback functor $\textbf{Vect}(X) \to \textbf{Vect}(X')$ is essentially surjective by Lemma \ref{finitedescent}. 

For the last statement, refer to the proof of \cite[Theorem 4.5]{SurfaceRes}: each vector bundle $\mathcal{F}_n$ on $Y$ fits in an exact sequence
\[0 \to \mathcal{F}_n \to (\mathcal{L}^{\otimes a_n})^{\oplus 2} \to \mathcal{L}^{\otimes b_n}|_{C_n}\to 0\]
where $a_n,b_n$ are positive integers, $\mathcal{L}$ is a line bundle and $C_n \subset Y$ is some Cartier divisor with $\mathcal{O}_{Y}(C_n) \cong \mathcal{L}^{\otimes 2a_n}$. By taking determinants, one sees that $\det(\mathcal{F}_n) \cong \mathcal{O}_Y$. Since $f^*\mathcal{E}_n \cong \mathcal{F}_n$ and $f$ is an isomorphism over $U$, the result follows.\end{proof}

\begin{lemma} \label{tensortracetrick} Let $F$ and $E$ be coherent sheaves on an algebraic space $X$. If $E$ is locally free, to show that $F$ admits a surjection by a locally free sheaf, it suffices to show that $F \otimes E$ admits a surjection from a locally free sheaf. 
\end{lemma}


\noindent The proof is straightforward so we omit it.

\begin{proof}[\textbf{Proof of Theorem \ref{ResSurface}}] If $X$ is not proper, then by \cite[Theorem 1.2.1]{CLO} there is a scheme-theoretically dense compactification $X \to \overline{X} \to \Spec A$. Moreover, since $X$ is catenary and Jacobson, \cite[Tag 0DS6]{stacks} implies $\dim(\overline{X})=2$. Thus, we may assume that $X$ is proper. We may also assume that $X$ is $S_1$ by passing to the $S_1$-ization $i: X' \to X$ (see \cite[Section 1]{NormalSurfaceRes}). Indeed, the algebraic space $X'$ satisfies the resolution property if and only if there is a $\GL_n$-torsor $W'$ on $X'$ with quasi-affine total space (see \cite[Theorem 1.1]{GrossRes}). However, any such torsor extends along $i$ because the sheaf of ideals, $I$, of the nilimmersion $i$, is supported in dimension $1$. Indeed, we may assume that $I$ is square zero and then the obstruction to deforming $W'$ lies in $\H^2(X, \text{Lie}(\text{Aut}(W')) \otimes I)=0$.

Let $M$ be a coherent sheaf, by Proposition \ref{freethick} there is a surjection $F \to M \to 0$ where $F$ is locally free on a thick open subset $U \subset X$. By lemma \ref{tensortracetrick}, it suffices to show that $F \otimes E$ admits a surjection from a vector bundle for \emph{some} rank $2$ locally free sheaf $E$. Moreover, by Proposition \ref{flat-glueing} there exists a coherent extension, $L$, of $\det(F|_U \otimes E|_U)^{\vee}$, and a class $\sigma \in \H^2(X, \Hom(F \otimes E, L))$ whose vanishing guarantees the existence of an extension
\[0 \to L \to N \to F \otimes E \to 0\]
where $N$ is $F_1$. Now we may write
\[\det(F|_U \otimes E|_U)^{\vee} \cong (\det(F|_U)^{\otimes 2})^{\vee} \otimes (\det(E)^{\otimes \text{rk}(F)})^{\vee}|_U\]
Let $L'$ be a fixed coherent extension of $(\det(F|_U)^{\otimes 2})^{\vee}$. Then, by \cite[Lemma 2.1]{SurfaceRes}, an isomorphism $\det(E)|_U \cong \mathcal{O}_X|_U$ induces an isomorphism on higher cohomology
\[\H^2(X, \Hom(F \otimes E,L)) \cong \H^2(X, \Hom(F,L \otimes E^{\vee})) \cong \H^2(X, \Hom(F, L') \otimes E^{\vee}))\]
Now Lemma \ref{1ample} implies there exists a rank $2$ vector bundle $\mathcal{E}_n$ with $\det(\mathcal{E}_n)|_U \cong \mathcal{O}_X|_U$ and
\[\H^2(X,\Hom(F, L') \otimes \mathcal{E}_n)=0\]
In conclusion, $F \otimes \mathcal{E}_n^{\vee}$ admits a surjection from a $F_1$ sheaf which is locally free on $U$. Thus, we may assume $M$ is a $F_1$ sheaf which is locally free on $U$. Next, use Lemma \ref{1ample} to find a rank $2$ vector bundle $\mathcal{E}_j$ with $\H^2(X, M^{\vee} \otimes \mathcal{E}_j)=0$, then for any $m>0$ apply \cite[Lemma 2.1]{SurfaceRes} to obtain isomorphisms
\[\H^2(X, \Hom(M, \mathcal{E}_j^{\oplus m})) \cong \H^2(X, M^{\vee} \otimes \mathcal{E}_j^{\oplus m}) \cong \H^2(X, M^{\vee} \otimes \mathcal{E}_j)^{\oplus m}=0\]
Now, Lemma \ref{Gross2} implies there is $m \gg 0$ and an extension
\[0 \to \mathcal{E}_j^{\oplus m} \to N \to M \to 0\]
where $N$ is locally free. \end{proof} 

\subsection{Many tame Artin gerbes satisfy the Resolution property}

In what follows we show that many tame Artin gerbes $\ms X$ satisfy the resolution property. Under the additional hypothesis that $\ms X$ is Deligne-Mumford with quasi-projective coarse moduli space this was proven in \cite{KreschVistoli} using a result of Gabber: $\Br(X)=\Br'(X)$ for all quasi-projective schemes $X$. We will show that if $\Br=\Br'$ for a class of algebraic stacks, then all tame Artin gerbes over such stacks admit (relatively) faithful vector bundles. We will also exercise control over the faithful vector bundle we produce on $\ms X$: their fibers will induce a power of the regular representation at every geometric point. As we will see in the next section, such vector bundles will be easy to glue over points where stacks are no longer gerbes. We begin with a definition of such vector bundles.

\begin{definition} \label{regularvb} Let $\ms X$ be a algebraic stack with finite diagonal. We say that a vector bundle $V$ of constant rank on $\ms X$ is \emph{regular} if for every geometric residual gerbe $\text{BG}_p$ the restriction of the vector bundle $V|_{\text{BG}_p}$ induces a power of the regular representation of the finite group scheme $G$. If $\ms X \to \ms Y$ is a morphism of algebraic stacks with finite diagonal, we say that a vector bundle on $\ms X$ is \emph{relatively regular} if it is a regular vector bundle after passing to a smooth schematic cover $Y \to \ms Y$. 
\end{definition}

\begin{proposition} \label{regulardescent} Let $f: \ms Y \to \ms X$ be a finite flat cover of constant degree between tame Artin gerbes which are finite-type over a field $k$. If $V$ is a regular vector bundle on $\ms Y$ then $f_*V$ is a regular vector bundle on $\ms X$. If $f: \ms Y \to \ms X$ is a finite flat cover of constant degree between tame Artin gerbes over an algebraic stack $\ms Z$ and $V$ on $\ms Y$ is relatively regular (with respect to $\ms Z$) then $f_*V$ is relatively regular (with respect to $\ms Z$).
\end{proposition}

\begin{proof} Let $\ms Y \to Y$ and $\ms X \to X$ denote the coarse moduli spaces and consider the factorization $\ms Y \to \ms X_Y \to \ms X$.  We will show that each of these morphisms is finite flat and that the pushforward of a regular vector bundle along each of them remains regular. It is not difficult to see that each of these maps is finite, so we explain why they are flat.

If $\ms Y \to \ms X_Y$ is finite flat then so is $\ms X_Y \to \ms X$. Indeed, $\ms Y \to \ms X$ is flat and $\ms X_Y \to \ms X$ being flat can be checked flat locally on the source. Flatness of $\ms Y \to \ms X_Y$ follows because it is a morphism between gerbes over the same algebraic space. It follows that $\ms X_Y \to \ms X$ is finite flat and therefore the induced map $Y \to X$ must also be finite flat, since $\ms X \to X$ is a smooth morphism. 

Since the morphism $\ms X_Y \to \ms X$ is affine, pushforward commutes with arbitrary base change on $\ms X$. Thus, we may assume that $X=\Spec \bar{k}$ and $Y$ is a finite $\bar{k}$-scheme. In this case $\ms X \cong \text{BG}_{\bar{k}}$ and $V$ is the pullback of a regular vector bundle on $\ms X$. Indeed, if $W$ is a regular vector bundle on $\ms X$ with the same rank as $V$ then $W|_{(\ms X_Y)_{red}} \cong V|_{(\ms X_Y)_{red}}$. But by Lemma \ref{vectorbundlesrigid}, it follows that $W|_{\ms X_Y} \cong V$. Thus, by the projection formula, the pushforward of $V$ is regular because it is a direct sum of copies of $W$. 

It remains to see why the pushforward of a regular vector bundle along $\ms Y \to \ms X_Y$ is regular. Since the morphism is finite the diagonal is a monomorphism and therefore, \'etale locally on $Y$, this looks like a morphism $\text{BH}_U \to \text{BG}_U$ where $\text{H}_U$ is a subgroup scheme of the finite flat group scheme $\text{G}_U$. Note that cohomology commutes with arbitrary base change on $\text{BG}_U$ since the map $\text{BH}_U \to \text{BG}_U$ is affine. As such, we need only check that pushing forward the regular representation along the map $\text{BH}_{\Spec L} \to \text{BG}_{\Spec L}$ yields a power of the regular representation. This follows since pushforward corresponds to taking the induced representation (see \cite[2.1]{Tame} and the preceding discussion) this is immediate (see also \cite[Chapter 1, 3.7.1]{Jantzen}). The statement about relatively regular vector bundles now follows by passing to a smooth schematic cover $Z \to \ms Z$. \end{proof}

\begin{lemma} \label{lem:outetale} Let $G$ be a locally well-split group scheme over a field $k$, then the outer automorphism group of $G$ is a finite \'etale group scheme over $k$. 
\end{lemma}

\begin{proof} By replacing $k$ with its algebraic closure, we may write $G=H \ltimes \Delta$ where $H$ is a tame constant group scheme and $\Delta$ is the connected component at the identity of $G$ (see \cite[Theorem 2.16]{Tame}). Then there exists a homomorphism of group schemes $\Delta \to \text{Aut}(G)$ which sends each section of $\Delta(T)$ to the corresponding inner automorphism of $G|_T$. The image of this group homomorphism is the connected component of \text{Aut}(G) by \cite[Lemma 2.19]{Tame}. It follows that $\text{Aut}(G)/\text{Inn}(G)=\text{Out}(G)$ is reduced over the algebraically closed field $k$ and so by \cite[Tag 047P]{stacks} the outer automorphism group is smooth. Moreover, finiteness of $G$ forces the finiteness of $\text{Aut}(G)$ and hence of $\text{Out}(G)$. This completes the proof. \end{proof}

Recall that if a gerbe $\ms X \to X$ becomes isomorphic to $\text{BG}_X$ locally on a site $X$, then the sheaf of isomorphisms of bands \[I=\text{Isom}(\text{Band}(G),\text{Band}(\ms X)) \to X\]
is a torsor under the outer automorphism group of $G$ (see e.g. \cite[3.2,3.3]{EHKV}).

\begin{lemma} \label{lem:bandconnected} Let $\ms X \to X$ denote a tame Artin gerbe over a noetherian and connected algebraic space $X$ over a field $k$. Then there is well-split group scheme $H \ltimes \Delta$ over $\Spec k$ with the property that after a \'etale base change $X' \to X$, we have $\text{Band}(\ms X)|_{X'} \cong \text{Band}(H \ltimes \Delta)|_{X'}$. 
\end{lemma}

\begin{proof} Note that a base change to the strict Henselization of any point $\Spec \mathcal{O}_p^{sh} \to X$ must trivialize the gerbe and yield an isomorphism 
\[\mathcal{X}|_{\Spec \mathcal{O}_p^{sh}} \cong \text{BG}_p|_{\Spec \mathcal{O}_p^{sh}}\]
for some finite flat linearly reductive group scheme $G_p \to \Spec \mathcal{O}_p^{sh}$. Now, by \cite[Theorem 2.16]{Tame}, we may deduce that $G_p$ is well-split after passing to a fppf cover $T \to \Spec \mathcal{O}_p^{sh}$:
\[G_p|_T \cong (H \ltimes \Delta)|_T\]
where $H \ltimes \Delta$ is a well-split group scheme over $k$. In particular, the sheaf 
\[I=\text{Isom}_{\text{Band}}(\text{Band}(\ms X)|_{\Spec \mathcal{O}_p^{sh}},\text{Band}(H \ltimes \Delta)|_{\Spec \mathcal{O}_p^{sh}}) \to \Spec \mathcal{O}_p^{sh}\] 
admits a $T$-point. However, this sheaf is a torsor under 
\[\text{Aut}_{\text{Band}}(\text{Band}(H \ltimes \Delta)|_{\Spec \mathcal{O}_p^{sh}})=\text{Out}(H \ltimes \Delta)|_{\Spec \mathcal{O}_p^{sh}}\]
and is therefore finite \'etale by Lemma \ref{lem:outetale}. Thus, $I$ admits a $\Spec \mathcal{O}_p^{sh}$-point as well. 

Thus, each $p \in X$ admits an \'etale neighborhood $U_p \to X$ so that $\text{Band}(\ms X)|_{U_p} \cong \text{Band}(G_p)|_{U_p}$ where $G_p$ is a well-split group scheme over $\Spec k$. It remains to show we may choose a single well-split group scheme that works for every $p \in X$. Since $X$ is connected it suffices to show that if $\eta$ specializes to $p$ in $X$ then $\text{Band}(G_p) \cong \text{Band}(G_{\eta})$ over $k$. 
However the sheaf of isomorphisms
\[I=\text{Isom}_{\text{Band}}(\text{Band}(G_{\eta}), \text{Band}(G_p)) \to \Spec k\]
is a psuedo-torsor under $\text{Aut}_{\text{Band}}(\text{Band}(G_p))=\text{Out}(G_p)$. Since both these bands are isomorphic over $\Spec \overline{k(\eta)}$, the scheme $I$ is nonempty, i.e. it is a torsor. However, since $\text{Out}(G_p)$ is a finite \'etale group scheme by Lemma \ref{lem:outetale}, the torsor $I$ becomes trivial after passing to a finite separable extension $k'/k$. The result now follows. \end{proof}

\begin{Theorem} \label{thm:tamegerbe} Let $\mathcal{S}$ denote a class of algebraic stacks finite-type and separated over a field $k$ with the following two properties 
\begin{enumerate}

\item If $X \in \mathcal{S}$ and $Y \to X$ is finite \'etale, then $Y \in \mathcal{S}$

\item If $X \in \mathcal{S}$ then $\Br(X)=\Br'(X)$

\end{enumerate}

\noindent Suppose $\ms X$ is a tame Artin gerbe over an object $X \in \mathcal{S}$. Then each connected component of $\ms X$ admits a vector bundle which is relatively regular over $X$. In particular, the morphism $\ms X \to X$ has the resolution property. \end{Theorem}

\begin{proof} We may assume that $k$ is separably closed by Proposition \ref{regulardescent}. We also assume that $\ms X$ is connected so by Lemma \ref{lem:bandconnected}, $\text{Band}(\ms X)$ is \'etale locally isomorphic to a well-split group scheme $G=H \ltimes \Delta$ defined over $k$. Thus, the sheaf of isomorphisms in the stack of bands over $X$:
\[I=\text{Isom}_{\text{Band}_X}(\text{Band}(\ms X), \text{Band}(H \ltimes \Delta)) \to X\]
is a $\text{Out}(G)|_X$-torsor. So if we replace $\ms X$ with $\ms X_I=\ms X \times_X I$ we may suppose that $\text{Band}(\ms X) \cong \text{Band}(H \ltimes \Delta)$. Indeed, $I \to X$ is finite \'etale by Lemma \ref{lem:outetale} and if $V$ is a relatively regular vector bundle on $\ms X_I$ (over $X$), the pushforward $(\text{pr}_1)_*V$ is a relatively regular vector bundle on $\ms X$ (over $X$) by Lemma \ref{regulardescent}.

Consider the stack of those equivalences which induce the identity on bands (see \cite[Chapter IV, Definition 2.3.1]{giraud}):
\[\ms H=\text{HOM}_{\text{id}}(\text{B}(H \ltimes \Delta), \ms X) \to X.\]
By \cite[Chapter IV, Theorem 2.3.2 (iii)]{giraud}, $\ms H$ is banded by $Z(H \ltimes \Delta) \subset H \ltimes \Delta$ and the natural evaluation map $\text{ev}: \ms H \to \ms X$ yields a map on bands induced by the inclusion $Z(H \ltimes \Delta) \subset H \ltimes \Delta$. However, because the base field is separably closed, the center of a well-split group scheme is diagonalizable. Thus we may write $\ms H \cong \ms R_1 \times_X \dots \times_X \ms R_n$ where each $\ms R_i$ is a $\mu_{n_i}$-gerbe over $X$. By hypothesis (1), $X$ remains in $\mathcal{S}$ and so by hypothesis (2) each $\ms R_i$ admits a $1$-twisted vector bundle, $V$, and therefore by taking an appropriate direct sum of copies of the vector bundles $V,V^{\otimes 2},...,V^{\otimes n_i}$ we may produce a regular vector bundle on $\ms R_i$. By pulling back these vector bundles to $\ms H$ via the projection maps and taking appropriate direct sums, it follows that $\ms H$ admits a relatively regular vector bundle. Moreover, fppf locally on $X$, the morphism $\text{ev}$ is of the form $\text{B}(Z(H \ltimes \Delta))|_X \to \text{B}(H \ltimes \Delta)|_X$, thus the pushforward of a relatively regular vector bundle along $\text{ev}$ remains relatively regular by Proposition \ref{regulardescent}. \end{proof}

\begin{proof}[\textbf{Proof of Theorem \ref{TAG}}] For (1), apply Theorem \ref{thm:tamegerbe} to the class $\ms S$ of quasi-projective schemes. Indeed, by Gabber's theorem (see \cite{dejongample}) we know $\Br(X)=\Br'(X)$ for every $X \in \ms S$. This shows every tame Artin gerbe over a quasi-projective scheme has the resolution property by Lemma \ref{resolutionpropertycomposes}. 

For (2), it suffices to show that (b) implies (a) since the converse is trivial. Let $\ms S$ (resp. $\ms S'$) denote the class of smooth (resp. normal) tame Artin stacks which have generically trivial stabilizers and a quasi-projective coarse moduli space. If (b) holds then any tame Artin gerbe morphism $\ms X \to X$ has the resolution property when $X \in \ms S$ (or $\ms S'$). Indeed, by (b) the stack $X$ has the resolution property and is therefore a global quotient stack. Thus, we may apply \cite[Theorem 4.4]{1RP} to obtain a finite flat cover $Z \to X$ where $Z$ is a quasi-projective scheme. Now Gabber's theorem (as in \cite{dejongample}) implies $\Br(Z)=\Br'(Z)$ and so Theorem \ref{BrauerMapTFAE} (5) forces $\Br(X)=\Br'(X)$. Thus, Theorem \ref{thm:tamegerbe} shows every tame Artin gerbe $\ms X \to X$ has the resolution property if $X \in \ms S$ (resp. $\ms S'$).

Now suppose that $\ms X$ is as in (a), then we may write $\ms X$ as a gerbe $\ms X \to \ms X^{\text{rig}}$ where $\ms X^{\text{rig}}$ is a smooth (resp. normal) tame Artin stack with generically trivial stabilizers by Theorem \ref{tamerigidification}. Since (b) holds, $\ms X^{\text{rig}}$ has the resolution property, and by the discussion above, the morphism $\ms X \to \ms X^{\text{rig}}$ has the resolution property. Thus, Lemma \ref{resolutionpropertycomposes} implies $\ms X$ satisfies the resolution property. \end{proof}

\subsection{Tame Artin Curves are Quotient Stacks}
In this section we show that tame Artin curves are always global quotient stacks and satisfy the resolution property. This has been known in the smooth case (in \cite[2.17]{EHKV} or \cite{BehrendNoohi} in the analytic or topological category) but there is no proof when the curve is allowed to be singular.

\begin{Theorem} \label{DMcurves} Let $\mathscr{X}$ be a tame Artin stack which is separated and finite-type over a field $k$ and suppose that all of the components of the coarse moduli space $X$ are of dimension $\leq 1$, then $\ms X$ is a global quotient stack.
\end{Theorem}

\begin{proof} 

It suffices to prove the result when $\ms X$ is connected. Next, we show it is enough to assume that $\mathscr{X}$ is reduced. Indeed, if $\mathscr{X}_{\text{red}}$ is a quotient stack, then it admits a faithful vector bundle $V$. We will deform along the nilpotent closed immersion $\mathscr{X}_{\text{red}} \to \mathscr{X}$ by factoring the extension into square zero extensions. Then the obstruction to deforming $V$ at each stage lies in $\H^2(\mathscr{X}_{\text{red}}, \text{End}(V) \otimes I)$ where $I$ is an ideal sheaf. Since $\mathscr{X}$ is tame, $\pi_*$ is exact and so the Leray spectral sequence for the coarse space map $\mathscr{X} \to X$:
\[\H^p(X, \R^q\pi_*( \text{End}(V) \otimes I)) \Rightarrow \H^{p+q} (\mathscr{X},  \text{End}(V) \otimes I)\]
degenerates. In particular we obtain isomorphisms 
\[\H^2(\mathscr{X},  \text{End}(V) \otimes I) \cong \H^2(X, \pi_*( \text{End}(V) \otimes I))\]
and the latter vanishes because $X$ is a scheme of dimension $1$.

Since $\ms X$ is reduced, the inertia stack $I_{\ms X} \to \ms X$ is flat on a dense open substack $\ms U$. Thus $\ms U$ is a gerbe and therefore by Theorem \ref{thm:tamegerbe} there exists a regular vector bundle on each component of $\ms U$. By taking appropriate direct sums of these regular vector bundles on the different components of $\ms U$, we may arrange for them to all have the same rank. Thus, we obtain a regular vector bundle $V$ on $\ms U$. Our aim is to extend a direct sum of $V$ over all the finitely many points of $\ms X \backslash \ms U$. By induction, it suffices to show that we can extend some sum $V^{\oplus m}$ over a single closed point of $\ms X \backslash \ms U$. Thus, we assume that $\ms X \backslash \ms U$ consists of one point, call it $p$ and let $R=\Spec \mathcal{O}_p^{sh} \to X$ denote the corresponding strictly local ring. Note that $R$ is one-dimensional since $X$ is a curve. Let $\ms U \to U$ denote coarse moduli, then there is a commutative diagram, all of whose squares are cartesian:

\[\begin{tikzcd}[row sep=scriptsize, column sep=scriptsize]
& \sX_{U'} \arrow[dl] \arrow[rr] \arrow[dd] & & \ms U \arrow[dl] \arrow[dd] \\
\sX_{R} \arrow[rr, crossing over] \arrow[dd] && \sX \\
& U' \arrow[dl] \arrow[rr] & & U \arrow[dl] \\
R \arrow[rr] & & X \arrow[from=uu, crossing over] \\
\end{tikzcd}\]

Moreover, by \cite[Theorem 3.2 (d)]{Tame} we can write $\ms X_R=[S/G]$ where $G$ is a locally well-split group scheme over $R$. Now, Theorem \ref{thm:tamegerbe} implies the pushforward of the structure sheaf along the map $S \to \ms X_R=[S/G]$ is a regular vector bundle on $\ms X_R$, call it $W$. By taking appropriate direct sums of $W$ and $V$, we may suppose that $W^{\oplus l}$ and $V^{\oplus m}$ have the same (constant) rank. By Theorem \ref{MVglueing}, to get a regular vector bundle on $\ms X$, it suffices to produce an isomorphism between $W^{\oplus l}$ and $V^{\oplus m}$ when restricted to $\ms X_{U'}$. Note that because $U'$ is a zero-dimensional noetherian scheme (because it is the complement of the closed point in $R$), it is the spectrum of Artinian local rings. Thus, since $V|_{\ms X_{U'}}^{\oplus m}$ and $W|_{\ms X_{U'}}^{\oplus l}$ are both regular vector bundles of the same constant rank, we may apply Lemma \ref{vectorbundlesrigid} to deduce that they are isomorphic over $U'$. The result follows. \end{proof}

\begin{corollary} \label{trivialdeterminant1} Let $\mathscr{X}$ be a tame Artin stack which is separated and finite-type over a field $k$ and suppose that all of the components of $X$ are of dimension $\leq 1$, then $\mathscr{X}$ admits a faithful vector bundle of constant rank with trivial determinant.
\end{corollary}

\begin{proof} By the previous theorem there exists a faithful vector bundle of constant rank on $\mathscr{X}$, call it $V$. Set $W=V \oplus \text{det}(V)^{\vee}$, it is still faithful since $V$ is but it also has trivial determinant since the determinant is additive in short exact sequences. \end{proof}

\begin{proof}[\textbf{Proof of Theorem \ref{T1} (1) and Theorem \ref{T2} (1)}] Let $\mathscr{X} \to X$ denote the coarse space map. By Theorem \ref{DMcurves} and Proposition \ref{resolutionpropertycoarsemoduli}, this morphism satisfies the resolution property. Since $X$ is $1$-dimensional, it is quasi-projective and thus satisfies the resolution property. By  Lemma \ref{resolutionpropertycomposes}, it follows that $\mathscr{X}$ has the resolution property. Also, any $\mu_n$-gerbe over a tame Artin curve is a tame Artin curve, and therefore has the resolution property as well. Theorem \ref{BrauerMapTFAE} immediately implies $\Br(\ms X)=\Br'(\ms X)$. \end{proof}

\section{The Surjectivity of the Brauer Map}

\subsection{Gluing Twisted Vector Bundles along flat neighborhoods}

We begin with a lemma which extends the following result from Gabber's thesis: if $X=U \cup V$ is a scheme covered by open affine sets $U,V$ with $U \cap V$ affine then $\Br(X)=\Br'(X)$.

\begin{lemma} \label{GlueingTwisted1} Suppose an algebraic stack $\mathscr{X}$ appears in a commutative diagram

\[\begin{tikzcd}
\Spec C \rar \dar & \mathscr{Y} \dar \\
\mathscr{X}' \rar & \mathscr{X} \\
\end{tikzcd}\]
such that every $\mu_n$-gerbe $\sG$ on $\sX$ satisfies the following two properties:

\begin{enumerate}
\item The natural functor
\[\mathbf{Vect}(\sG) \to \mathbf{Vect}(\sG_{\mathscr{X}'}) \times_{\mathbf{Vect}(\sG_{\spec(C)})} \mathbf{Vect}(\sG_{\mathscr{Y}})\]
is an equivalence of categories.

\item There exists twisted vector bundles of the same nonzero constant rank $V$ and $W$ on $\sG_{\sX'}$ and $\sG_{\sY}$ respectively which satisfy $V|_{\sG_{\Spec C}} ^{\otimes n} \cong W|_{\sG_{\Spec C}} ^{\otimes n} \cong \mathcal{O}_{\sG_{\Spec C}}^{\oplus m}$ for some positive natural number $m$.
\end{enumerate}

Then every $\mu_n$-gerbe on $\sX$ admits a twisted vector bundle of nonzero constant rank, and in particular $\Br(\sX)=\Br'(\sX)$.

\end{lemma}

\begin{proof}
Let $\sG$ denote a $\mu_n$-gerbe on $\mathscr{X}$, it suffices to exhibit a nonzero twisted vector bundle of constant rank on $\sG$. Let $V$ and $W$ be twisted vector bundles as in (2). If we show that some direct sum of $V$ agrees with a direct sum of $W$ on $\sG_{\Spec C}$, then (1) allows us to conclude that there is a twisted vector bundle on $\sG$.

Observe that on $\sG_{\Spec C}$ both the untwisted modules $W|_{\sG_{\Spec C}} \otimes V|_{\sG_{\Spec C}}^{\vee}$, $V|_{\sG_{\Spec C}} \otimes V|_{\sG_{\Spec C}}^{\vee}$ correspond to (the pullback of) projective modules on $\Spec C$. Moreover, their $n$th-tensor powers are both free by (2). It follows from K-theory as in \cite[Proposition 3.1.4.3 (iv)]{lieblich2008twisted} that there is a $N$ such that $(W|_{\sG_{\Spec C}} \otimes V|_{\sG_{\Spec C}}^{\vee})^{\oplus N} \cong (V|_{\sG_{\Spec C}} \otimes V|_{\sG_{\Spec C}}^{\vee} )^{\oplus N} \cong F$ where $F$ is a free $C$-module of rank $M>0$. Also note that there is an isomorphism due to associativity
\[V|_{\sG_{\Spec C}} \otimes (W|_{\sG_{\Spec C}} \otimes V|_{\sG_{\Spec C}}^{\vee}) \cong W|_{\sG_{\Spec C}} \otimes (V|_{\sG_{\Spec C}} \otimes V|_{\sG_{\Spec C}}^{\vee})\]
Thus we obtain

\begin{align*}
V|_{\sG_{\Spec C}}^{\oplus M}
& \cong
V|_{\sG_{\Spec C}} \otimes F \\
& \cong
V|_{\sG_{\Spec C}} \otimes (W|_{\sG_{\Spec C}} \otimes V|_{\sG_{\Spec C}}^{\vee})^{\oplus N}\\
& \cong
[V|_{\sG_{\Spec C}} \otimes (W|_{\sG_{\Spec C}} \otimes V|_{\sG_{\Spec C}}^{\vee})]^{\oplus N} \\
& \cong
 [W|_{\sG_{\Spec C}} \otimes (V|_{\sG_{\Spec C}} \otimes V|_{\sG_{\Spec C}}^{\vee})]^{\oplus N} \\
& \cong
W|_{\sG_{\Spec C}} \otimes (V|_{\sG_{\Spec C}} \otimes V|_{\sG_{\Spec C}}^{\vee})^{\oplus N} \\
& \cong
W|_{\sG_{\Spec C}} \otimes F \\
& \cong
W|_{\sG_{\Spec C}}^{\oplus M} \\
\end{align*}

The result follows. \end{proof}

\subsection{A Stacky generalization of Schr\"{o}er's Theorem}

The purpose of this section is generalize Schr\"oer's Theorem: there are enough Azumaya algebras on normal surfaces. Indeed, this is equivalent to the following: if $f:\ms X \to X$ is a $\mu_n$-gerbe over a normal surface, then $f$ has the resolution property. On the other hand, Theorem \ref{BrOrbisurface} (see below) shows that any tame Artin stack which is a gerbe over an $R_1$ orbisurface must have the resolution property. An immediate corollary is that \emph{all} normal tame Artin surfaces have the resolution property.

\begin{lemma} \label{lifttrivial} Let $\mathscr{X}_{\Spec A} \to \Spec A$ be a coarse moduli space of a tame Artin stack where $A$ is noetherian and $I$-adic for some ideal $I$. Suppose moreover that there exists a vector bundle, $V_0$, on $\mathscr{X}_{\Spec A/I}$ with trivial determinant. Then $\mathscr{X}$ admits a vector bundle $V$ with trivial determinant lifting $V_0$.
\end{lemma}

\begin{proof}
To obtain such a vector bundle on $\mathscr{X}_{\Spec A} \to \Spec A$ we use deformation theory to obtain a compatible collection of vector bundles over $\{\mathscr{X}_{\Spec A/I^n}\}$ and then algebraize using a generalization of Grothendieck's existence theorem.
We deform the vector bundle $V_0$ along the sequence of thickenings
\[\mathscr{X}_{\Spec A_0} \to \mathscr{X}_{\Spec A/I^2} \to \mathscr{X}_{\Spec A/ I^3} \to ...\]
At each stage a deformation $(V_i, \phi)$ on $\mathscr{X}_{\Spec A/I^{i+1}}$ exists and is unique up to isomorphism by Proposition \ref{deformationobstruction}, the Leray spectral sequence for $\pi: \mathscr{X}_{\Spec A_0} \to \Spec A_0$, and tameness. Indeed the deformation and obstruction spaces are
\[\H^i(\mathscr{X}_{\Spec A_0}, \End(V_{i-1}) \otimes I^{i+1}/I^{i+2}) \cong \H^i(\Spec A_0, \pi_* (\End(V_{i-1}) \otimes I^{i+1}/I^{i+2}))\]
for $i=1, 2$ respectively and these groups vanish because an affine scheme has no quasi-coherent cohomology. Note that a deformation of a vector bundle $V_i$ also induces a deformation of its determinant, but since the deformation spaces are trivial and $\text{det} (V_0) \cong \mathcal{O}_{\mathscr{X}_{\Spec A_0}}$ this implies $\det(V_1) \cong \mathcal{O}_{\mathscr{X}_{\Spec A/I^2}}$. By induction it follows that the $V_i$ all have trivial determinants.
Thus we have a compatible system $\{V_i\}$ of vector bundles all of whose determinants are trivial. By a generalization of Grothendieck's Existence theorem (see \cite[Theorem 1.4]{olssonproper}) there exists a vector bundle on $\mathscr{X}_{\Spec A}$ restricting to $\{V_i\}$ with trivial determinant. \end{proof}

\begin{Theorem} \label{BrOrbisurface} Let $\mathscr{X}$ be a $2$-dimensional, tame Artin stack which is finite-type and separated over a field. Suppose that $\ms X$ is a gerbe over a stack $\ms X'$ with generically trivial stabilizers. If $\ms X$ is regular in codimension $1$ ($R_1$), then it has the resolution property.
\end{Theorem}

\begin{proof} Since the gerbe morphism $\ms X \to \ms X'$ is smooth and codimension-preserving, it follows that $\ms X'$ is also $R_1$. Thus, to prove the theorem it suffices to show 
\begin{enumerate} \item $\ms X'$ is a global quotient stack, and

\item There are enough Azumaya algebras on $\ms X'$ i.e. $\Br(\ms X')=\Br'(\ms X')$. \end{enumerate}
Indeed, Theorem \ref{ResSurface}, Proposition \ref{resolutionpropertycoarsemoduli} and Lemma \ref{resolutionpropertycomposes} shows that (1) implies $\ms X'$ has the resolution property. Thus, to show the theorem it suffices to show the gerbe morphism $\ms X \to \ms X'$ has the resolution property. However, since the class of $R_1$ tame Artin surfaces with generically trivial stabilizers is stable under finite \'etale covers, Theorem \ref{thm:tamegerbe} implies it is enough to show $\Br(\ms X')=\Br'(\ms X')$. 

Let $\sG$ be a $\mu_n$-gerbe on $\sX'$ and denote by $\sX' \to X$ the coarse space map. It suffices to find a twisted vector bundle of nonzero rank on $\sG$. Then, Theorem \ref{BrauerMapTFAE} implies (2) holds. To show (1), we will exhibit a faithful vector bundle on $\ms X'$. Since the constructions of these two vector bundles are so similar, we explain them in parallel.

Let $\Spec A \subset \sX'$ be a dense smooth affine open substack and denote by $\sC$ the complement with its reduced induced structure and $\sC \to C$ the associated coarse space map. Since the (finite) singular locus of $\sX'$ lies on $\sC$, the singular points of $\ms X'$ all have an image in $C$, we denote them by $q_1,...,q_n$. 

Let $\sG_{\sC} \to \sC$ be the restricted $\mu_n$-gerbe, we begin by showing there exists a twisted vector bundle on $\sG_{\sC}$. Since $\sG_{\sC}$ is $1$-dimensional, tame, finite-type and separated over a field it has the resolution property by Theorem \ref{T1} (1). Now, Theorem \ref{BrauerMapTFAE} implies $\sG_{\sC}$ admits a twisted vector bundle of constant (nonzero) rank $V_{\sC}$. By replacing $V_{\sC}$ with direct sums of itself, we may suppose that $\det(V|_{\sC})$ is the pullback of a line bundle on $C$ (see \cite{RydhMO}). Moreover, by tensoring $V_{\sC}$ with a large power of an ample line bundle on $C$ we may also suppose that $\det(V|_{\sC})$ is very ample on $C$. Thus, by taking the complement of an appropriate section of $\det(V|_{\sC})$, we may find a open affine neighborhood $\Spec B_1 \subset C$ containing the $q_1,...,q_n$ and which also trivializes $\det(V|_{\sC})$. By \cite[Tag 04CE]{stacks}, there exists a open subspace $W_1 \subset X$ such that $W_1 \cap C=\Spec B_1$ so that the map $\Spec B_1 \to W_1$ is a closed immersion. Similarly, we may also find a $\Spec B_2 \subset C$, an open subspace $W_2 \subset X$ (both containing the $q_1,...,q_n$) so that $\Spec B_1 \cup \Spec B_2 = C$ and $W_2 \cap C =\Spec B_2$. 

Next, we will show that $\sG|_{W_1 \cup \Spec A}$ admits a twisted vector bundle of nonzero rank and that each of the open substacks, $\ms X'|_{W_i \cup \Spec A}$, admits a faithful vector bundle. Then the pushforward of these vector bundles to all of $\sG$ (or $\ms X'$) is locally free because a reflexive module on a regular local ring of dimension $2$ is free. Moreover, the pushforward of a twisted vector bundle is twisted by Lemma \ref{twistedopenclosed}. Lastly, the direct sum of the resulting vector bundles on $\ms X'$ is faithful since it is faithful over each member of the open cover $\ms X'=\ms X'_{W_1 \cup \Spec A} \cup \ms X'_{W_2 \cup \Spec A}$. 

First, use Theorem \ref{flatnhbdsexist} to obtain flat neighborhoods $\Spec R_i \to W_i \cup \Spec A$ of $\Spec B_i$ which are affine, adic and $2$-dimensional (since it surjects onto height $2$ points of $X$). Then, applying Lemma \ref{lifttrivial} to the coarse space morphism $\sG_{\Spec R_1} \to \Spec R_1$ and the twisted vector bundle $V_{\sC}|_{\Spec R_1}$, we obtain a \emph{twisted} vector bundle (see Proposition \ref{twistedopenclosed} for twistedness) $V_{\sG_{\Spec R_1}}$ on $\sG_{\Spec R_1}$ with trivial determinant. Similarly, we obtain faithful vector bundles, $D_i$, of constant rank on $\ms X'|_{\Spec R_i}$ with trivial determinant using Lemma \ref{lifttrivial} and Corollary \ref{trivialdeterminant1}.

Observe that we obtain the following diagram all of whose squares are cartesian and whose bottom (and top) faces are flat Mayer-Vietoris squares.

\[\begin{tikzcd}[row sep=scriptsize, column sep=scriptsize]
& \Spec A \times_X \Spec R_1 \arrow[dl] \arrow[rr] \arrow[dd] & & \sX_{\Spec R_i}' \arrow[dl] \arrow[dd] \\
\Spec A \arrow[rr, crossing over] \arrow[dd] && \sX_{\Spec A \cup W_i}' \\
& \Spec A \times_X \Spec R_i \arrow[dl] \arrow[rr] & & \Spec R_i \arrow[dl] \\
\Spec A \arrow[rr] & & \Spec A \cup W_i \arrow[from=uu, crossing over] \\
\end{tikzcd}\]

Note that the restriction $D_i|_{\Spec A \times_X \Spec R_i}$ is a free module by Serre's splitting theorem (\cite[Theorem 1]{SerreSplitting}) since $\det(D_i)|_{\Spec R_i}$ is free and $\Spec A \times_X \Spec R_i=\Spec R_i \backslash \Spec B_i$ is a $1$-dimensional affine scheme by Theorem \ref{flatnhbdsexist}. Thus, $D_i \in \textbf{Vect}(\ms X'_{\Spec R_i})$ extends to a faithful vector bundle on $\ms X'_{\Spec A \cup W_i}$ by Theorem \ref{MVglueing}. The statement (1) follows. 

We may find a constant nonzero rank twisted vector bundle on $\sG_{\Spec A}$ whose $n$th tensor power is free. Indeed, a nonzero twisted vector bundle, $J$, exists on $\sG_{\Spec A}$ by \cite[Theorem 1, Chapter 2]{Gabbersthesis}. Then, by \cite[Corollary 3.1.4.4]{lieblich2008twisted} we may tensor $J$ by an appropriate projective module so that $J^{\otimes n}$ is free. Indeed, to use the cited corollary set $P=J^{\otimes n}$ (since it is untwisted, it may be viewed as a locally free $A$-module) and replace $J$ with $J \otimes \bar{P} \otimes F_0$ where $\bar{P}$ and $F_0$ are as in the corollary. The $n$th tensor power of this vector bundle is free and $J$ itself is twisted since it is a twisted vector bundle tensored with untwisted modules.

By taking appropriate direct sums of $J$ and $V_{\sG_{\Spec R_i}}$, we may assume they have the same rank without changing the fact that $V_{\sG_{\Spec R_i}}$ has trivial determinant. We will glue these vector bundles along the flat Mayer-Vietoris square above. Indeed, to apply Lemma \ref{GlueingTwisted1}, it remains to see why $V_{\sG_{\Spec R_i \times_X \Spec A}}^{\otimes n}$ is trivial on $\Spec R_i \times_X \Spec A$. However, $\Spec R_i \times_X \Spec A$ is a noetherian $1$-dimensional affine scheme. Thus, Serre's splitting theorem (see \cite[Theorem 1]{SerreSplitting}), and the triviality of $\det(V|_{\sG_{\Spec R_i}})$ implies $ V_{\sG_{\Spec R_i \times_X \Spec A}}^{\otimes n}$ is isomorphic to a free module. The result follows. \end{proof}

\begin{proof}[\textbf{Proof of Theorem \ref{T1}(2) and Theorem \ref{T2}(2)}] If $\ms X$ is as in the statement of the Theorem \ref{T1} (2), then we may write it as a gerbe $\ms X \to \ms X'$ over a tame Artin surface $\ms X'$ with generically trivial stabilizers by Theorem \ref{tamerigidification}. The result now follows from Theorem \ref{BrOrbisurface}. Moreover, every $\mu_n$-gerbe over a $R_1$ tame Artin surface which is a gerbe over a stack with generically trivial stabilizers is, itself, also a gerbe over a stack with generically trivial stabilizers. Thus, Theorem \ref{BrauerMapTFAE} and Theorem \ref{BrOrbisurface} implies $\Br(\ms X)=\Br'(\ms X)$, as desired. \end{proof}

\label{sec:bibliography}

\bibliography{mybib}{}
\bibliographystyle{plain}

\end{document}